\documentclass[11pt,a4paper]{article}
\usepackage{fullpage}
\usepackage[margin=2.5cm]{geometry}
\usepackage{amsmath,amssymb,amsfonts,amsthm}
\usepackage[bookmarksnumbered,unicode]{hyperref}
\usepackage[capitalise]{cleveref}
\newtheorem{theorem}{Theorem}[section]
\newtheorem{lemma}[theorem]{Lemma}
\newtheorem{corollary}[theorem]{Corollary}
\newtheorem*{localclaim}{Claim}
\DeclareMathOperator{\tail}{tail}
\newcommand{\Q}{\mathbb Q}

\newcommand{\R}{\mathbb R}
\newcommand{\Z}{\mathbb Z}
\newcommand{\E}{\mathbb E}
\newcommand{\D}{\mathcal D}
\newcommand{\A}{\mathcal A}

\newcommand{\F}{\mathcal F}

\newcommand{\ex}{\mathrm{ex}}
\newcommand{\co}{\mathrm{co}}
\title{Prescribed number fields in single degree Tur\'an spectra}
\author{Yang Xiao\thanks{\raggedright\hangindent=1.8em\relax Department of Math and Computer Science, Biola University, La Mirada, CA 90639, USA.\newline Email: \texttt{yang.xiao@Biola.edu}.}
\and Haotian Yang\thanks{\raggedright\hangindent=1.8em\relax Department of Mathematics, California Institute of Technology, Pasadena, CA 91125, USA.\newline Email: \texttt{hyang3@caltech.edu}.}}
\date{}
\hypersetup{pdftitle={Prescribed number fields in single degree Turán spectra},pdfauthor={Yang Xiao, Haotian Yang}}

\begin{document}
\maketitle
\begin{abstract}
The \(\ell\)-degree Tur\'an density of a \(k\)-graph \(F\) is the limiting
maximum of \(\delta_\ell(H)/\binom{n-\ell}{k-\ell}\) over \(n\)-vertex \(k\)-graphs
\(H\) avoiding \(F\), where \(\delta_\ell(H)\) is the minimum, over all
\(\ell\)-sets \(S\), of the number of edges containing \(S\). We show that for every \(k\ge3\), \(1\le\ell<k\)
and number field \(K\subseteq\R\), there are such densities \(\alpha\)
arbitrarily close to \(1\) from below with \(\Q(\alpha)=K\). In particular,
for every \(k\ge3\), each positive integer occurs as the algebraic degree
of the codegree Tur\'an density of a single \(k\)-graph.
\end{abstract}

\section{Introduction}\label{ch:overview}

Fix integers \(k\ge3\) and \(1\le\ell<k\). For an \(\ell\)-set \(S\) in a \(k\)-uniform
hypergraph (or \(k\)-graph) \(H\), let \(d_H(S)\) be the number of edges
containing \(S\). The minimum \(\ell\)-degree of \(H\) is
\(\delta_\ell(H)=\min_{S\in\binom{V(H)}{\ell}}d_H(S)\).
Let \(\mathcal F\) be a nonempty family of \(k\)-graphs, each with at least one edge.
A \(k\)-graph is \(\mathcal F\)-free if it contains no member of \(\mathcal F\)
as a subhypergraph. The \emph{\(\ell\)-degree Tur\'an number}
\(\ex_\ell(n,\mathcal F)\) is the maximum of \(\delta_\ell(H)\) over all
\(\mathcal F\)-free \(k\)-graphs \(H\) on \(n\) vertices. The corresponding
\emph{\(\ell\)-degree Tur\'an density} is
\[
  \pi_\ell(\mathcal F)
  =\lim_{n\to\infty}\frac{\ex_\ell(n,\mathcal F)}{\binom{n-\ell}{k-\ell}}.
\]
Lo and Markstr\"om~\cite{LM} proved that this limit exists. When
\(\mathcal F=\{F\}\), we write \(\pi_\ell(F)\) for \(\pi_\ell(\mathcal F)\).
A vertex-deletion argument gives \(\pi_1(F)=\pi(F)\), the ordinary Tur\'an density.
For \(\ell=k-1\) we use the codegree notation \(\delta_{\co}=\delta_{k-1}\)
and \(\pi_{\co}=\pi_{k-1}\).

Determining the exact Tur\'an density of a given hypergraph is generally
difficult. Even the ordinary Tur\'an density of the complete \(k\)-graph
\(K_t^{(k)}\) on \(t\) vertices is unknown for every \(t>k\ge3\). This has led
to extensive work on the complementary question of which real numbers in
\([0,1)\) can occur as Tur\'an densities.

When the forbidden family may be infinite, Pikhurko~\cite{P} showed that
the set of ordinary Tur\'an densities has the cardinality of the continuum.
Liu and Pikhurko~\cite{LP} showed that it contains an interval \([c_k,1)\)
for some \(c_k<1\). For \(2\le\ell<k\), Li and Liu~\cite{LL} proved that
every \(\alpha\in[0,1)\) is an \(\ell\)-degree density. When \(\ell>2\),
Lin and Zhou~\cite{LZ} constructed a family of density \(\alpha\) whose
finite subfamilies all have strictly larger density. Gao, Pikhurko, Rong
and Sun~\cite{GPRS} proved the same for codegree in \(3\)-graphs.

For finite forbidden families, Baber and Talbot~\cite{BaberTalbot2012}
constructed irrational ordinary Tur\'an densities in uniformity three, and
Pikhurko~\cite{P} did so for every \(k\ge3\). Answering a question of
Grosu~\cite{Grosu}, Liu and Pikhurko~\cite{LPDegree} showed that such
densities have unbounded algebraic degree for every \(k\ge3\). Mubayi and Zhao~\cite{MZ}
showed that the codegree densities of finite families are dense in
\([0,1)\), and Lo and Markstr\"om~\cite{LM} extended this to every
\(2\le\ell<k\). Gao, Pikhurko, Rong and Sun~\cite{GPRSfinite} realized
every rational number in \([0,1)\) as a codegree density. Tang and
Zhou~\cite{TZ} proved the corresponding result for every \(k/2<\ell<k\).

We now turn to single forbidden \(k\)-graphs. Write \(\Pi^k_\ell\) for the set
of densities \(\pi_\ell(F)\) of finite \(k\)-graphs \(F\) with at least
one edge, and \(\Pi^k=\Pi^k_1\) for the ordinary spectrum. Conlon and
Sch\"ulke~\cite{CSfirst,CS} found infinitely many accumulation points of
\(\Pi^k\). Chen and Liu~\cite{ChenLiuClosure} showed that
\(\overline{\Pi^k}=\{0\}\cup[k!/k^k,1]\) for all sufficiently large \(k\). Examples of irrational ordinary densities were given by Yan and Peng~\cite{YanPeng}
and Wu~\cite{Wu}. Kam\v{c}ev, Letzter and Pokrovskiy~\cite{KLP} proved
that the tight \(3\)-uniform cycle of length \(t\) has density \(2\sqrt3-3\)
for every sufficiently large \(t\) not divisible by three.

Using tree suspensions, Ai, Ding, Liu and Yang~\cite{ADLY} showed that
\(\Pi^k_\ell\) contains algebraic numbers of arbitrarily large degree for
\(\ell=1\) and for \(k/2\le\ell\le k-2\). They asked whether the
algebraic degrees in the codegree spectrum are also unbounded
\cite[Problem~4.2]{ADLY}. Our main theorem answers their question and, for each
\(1\le\ell<k\), prescribes the field generated by a density in \(\Pi^k_\ell\).

\begin{theorem}\label{thm:main}
Let \(k\ge3\), \(1\le\ell<k\), and \(K\subset\R\) be a number field. For every \(\varepsilon>0\) there is a single finite \(k\)-graph \(F\) with \(1-\varepsilon<\pi_\ell(F)<1\) and \(\Q(\pi_\ell(F))=K\).
\end{theorem}

Here \(K\) is taken with its given embedding in \(\R\). For instance,
taking \(K=\Q(\sqrt[d]{2})\) gives a density of degree exactly \(d\)
for every positive integer \(d\), since \(X^d-2\) is irreducible over \(\Q\).

The rest of the paper is organized as follows.
In Section~\ref{sec:preliminaries}, we prove a product bound for passing
from a finite forbidden family to a single hypergraph, collect the
algebraic and colouring lemmas, and extend the tree suspension to finite families.
In Section~\ref{sec:ordinary-cases}, we prove the ordinary-density case
by adjoining independent classes to finite patterns and realizing the
resulting densities with single forbidden hypergraphs.
In Section~\ref{ch:cg}, we construct finite families of \(3\)-graphs
whose codegree density generates the prescribed field, lift the construction to
higher uniformities, and pass to a single forbidden hypergraph.
In Section~\ref{sec:intermediate-cases}, we treat the intermediate
degrees other than \((k,\ell)=(4,2)\) by expanding edges when
\(k-\ell<\ell\) and prescribing links when \(k-\ell\ge\ell\).
In Section~\ref{ch:quartic}, we use the lines of an affine space over
\(\mathbb F_4\) to settle the case \((k,\ell)=(4,2)\) and complete the
proof of Theorem~\ref{thm:main}.

\section{Preliminaries}\label{sec:preliminaries}

We first record the product argument that we shall use for the upper bounds. For a set \(S\subseteq V(H)\), let \(L_H(S)\) be the hypergraph on \(V(H)\setminus S\) with edges \(e\setminus S\) over all \(e\in E(H)\) containing \(S\). If \(H\) is an \(r\)-graph and \(|S|=r-1\), we write \(N_H(S)=\{v:S\cup\{v\}\in E(H)\}\) for its neighbourhood. For a single vertex we write \(L_H(v)\). In a \(3\)-graph we write \(N_H(x,y)\) for the neighbourhood of the pair \(\{x,y\}\). In a partitioned hypergraph we call a set \emph{pure} if all its vertices lie in one part, and \emph{mixed} otherwise.
A homomorphism of \(r\)-graphs sends each edge bijectively onto an edge. Write \(B(m)\) for the blow-up of \(B\) with \(m\) clones at every vertex. The categorical product of a family of \(r\)-graphs has the Cartesian product as its vertex set. An \(r\)-set is an edge precisely when its projection is an edge in every factor.

\begin{lemma}\label{lem:product-upper}
If \(P\) maps homomorphically to every member of a finite family \(\mathcal B\), then \(\pi_\ell(P)\le\pi_\ell(\mathcal B)\). In particular one can use \(P=\prod_{B\in\mathcal B} B\).
\end{lemma}
\begin{proof}
For each \(B\in\mathcal B\), fix a homomorphism \(P\to B\).
Sending vertices in the same fibre to distinct clones embeds \(P\) in \(B(|V(P)|)\).
Thus every \(P\)-free hypergraph avoids the whole family
\(\{B(|V(P)|):B\in\mathcal B\}\). Blow-up invariance for finite families
gives the inequality~\cite[Proposition~1.6]{LM}. For a product, the
homomorphisms are the coordinate projections.
\end{proof}
For related uses of this argument see~\cite{PS,DLLWY,ADLYzero,LiLiu}. In each application below the lower bound comes from a construction that avoids \(P\) and attains the density of the finite family.

Throughout, \(K\) is a number field with its given real embedding. Let
\(K_\infty=K\otimes_\Q\R\), with one coordinate for each real embedding
and one for each pair of nonreal conjugate embeddings. We shall need to
choose a generator with prescribed open conditions on its conjugates and
to recognize transformations that preserve the field it generates.

\begin{lemma}\label{lem:choose-generator}
Each nonempty open subset of \(K_\infty\) contains \(w\in K\) with distinct squares of its conjugates; in particular \(\Q(w)=K\).
\end{lemma}
\begin{proof}
The rational span of a basis of \(K\) is dense in \(K_\infty\). For each pair of distinct embeddings \(\sigma,\tau\), the equalities \(\sigma(w)=\tau(w)\) and \(\sigma(w)=-\tau(w)\) each cut out a proper real linear subspace. We may thus choose \(w\) outside their finite union in the given open set. Its conjugates are distinct, so \(w\) generates \(K\).
\end{proof}
\begin{lemma}\label{lem:field-transfers}
Let \(\Q(a)=K\). If \(M\in\Q(X)\) is a nonconstant fractional-linear
function defined at \(a\), then \(\Q(M(a))=K\). For each nonconstant \(R\in\Q(X)\), also \(\Q(R(q+a))=K\) for all but finitely many integers \(q\).
\end{lemma}
\begin{proof}
The inverse of \(M\) expresses \(a\) as a rational function of \(M(a)\). For the second assertion, let \(a_i,a_j\) be distinct conjugates. The rational functions \(R(X+a_i)\) and \(R(X+a_j)\) are distinct, since a nonconstant rational function has no nonzero period. Only finitely many integers \(q\) are poles or give a coincidence for some pair of conjugates~\cite[Lemma~3.1]{ADLY}.
\end{proof}

A weak \(p\)-colouring is a partition into \(p\) independent classes.
For a class \(\D\), a finite family of witnesses consists of hypergraphs
outside \(\D\). The family is \emph{exact} if its Tur\'an density for the
degree in question equals the asymptotic extremal density of \(\D\).
When \(\D\) is closed under taking subhypergraphs, no member of \(\D\)
contains a witness as a subhypergraph. A \emph{positive coloured prohibition}
prescribes vertex colours and a set of required edges. It occurs whenever
there is a colour-preserving copy of this subhypergraph.

To find witnesses we first colour a bounded set of vertices, the
\emph{anchor}. A \emph{certificate} for a colour restriction is a bounded
set of vertices whose incidences, together with the coloured anchor,
already force that restriction. A vertex is \emph{pinned} if such a
certificate forces its colour. The bounds here are independent of
the order of the host.
\begin{lemma}[Bounded colour certificates]\label{lem:bounded-colour}
Fix a bounded anchor, finitely many colours and positive coloured
prohibitions of bounded order. Fix also a colouring of the anchor.
Suppose that pins and auxiliary restrictions have bounded certificates
and that any two assignments satisfying the auxiliary restrictions on
the same vertex set differ on at most \(B\) unpinned vertices. If the
host has no representation extending the anchor colouring, a bounded
subhypergraph already has no such representation.
\end{lemma}
\begin{proof}
If the anchor itself fails a restriction, it is already a witness.
First suppose that the auxiliary restrictions are inconsistent. Take a
minimal inconsistent set of vertices outside the fixed anchor. Remove
one vertex \(v\) and choose a colouring of the remainder satisfying the
auxiliary restrictions. For each possible colour of \(v\), extend this
colouring to the whole set and keep the colour of \(v\) fixed. At every
stage record a violated auxiliary restriction and its certificate. Branch
by changing one as yet unchanged, unpinned vertex in its support, other
than \(v\), to each of its other colours. Stop after \(B\) changes and
record a violation also at each leaf.

Now suppose that the auxiliary restrictions are consistent. Start from a
colouring of the host satisfying those restrictions and all pins. When
a colouring satisfies the auxiliary restrictions, it violates an
original coloured prohibition since there is no representation. Record
one such prohibition. If a change has instead violated an auxiliary
restriction, record that restriction and its certificate. Use the same
branching rule up to depth \(B\), fixing only the anchor and the pins.

In either case let \(U\) consist of the vertices encountered in the
search, including \(v\) in the first case. Retain \(U\), the anchor and
the certificates for all pins and auxiliary restrictions on \(U\).
The depth, the number of colours and the sizes of the supports are
bounded, so this uses boundedly many vertices. Suppose a valid colouring existed there.
Follow the branch that agrees with this colouring on \(v\) in the first
case. At any recorded violation there must be an as yet unchanged
vertex whose colour differs from the proposed valid colouring. Change
it to that colour. On \(U\) in the second case, and on \(U\setminus\{v\}\) in the first,
both the proposed colouring and the baseline satisfy the auxiliary
restrictions. They differ on at most \(B\) vertices there. The violation recorded
at depth \(B\) therefore gives a contradiction.

If the host has no representation, taking the union over the boundedly
many anchor colourings gives a
bounded witness. In particular, for a class closed under taking
subhypergraphs, suppose the argument applies to all sufficiently large
nonmembers above a degree threshold below the class density. The
nonmembers up to the resulting order bound then form an exact finite
forbidden family.
\end{proof}

We shall also need to choose an ordered clique in each factor of a product so that the vertex in any fixed position is uniformly distributed.
\begin{lemma}[Uniform clique positions]\label{lem:uniform-cliques}
Let an \(r\)-graph of order \(m\) have at most \(\mu m^{r-1}\) missing edges at each vertex. If \(L\mid m\) and \(\mu L^{r-1}\) is small enough in terms of \(r\), it admits a uniform-marginal law on ordered \(L\)-cliques.
\end{lemma}
\begin{proof}
Split a random permutation of the \(m\) vertices into consecutive blocks of size \(L\). For each nonedge and each ordered choice of positions in a block, the corresponding bad canonical event has probability \(1/(m)_r\). Two such events can be joined whenever they share a vertex or a position. This gives a lopsided dependency supergraph~\cite{LuSzekely} of maximum degree \(O_r(\mu m^rL^{r-1})\), by counting assignments through a given vertex or position. The local lemma gives a partition in which every block is a clique. Choose a block uniformly and then order its vertices uniformly. Each position has the desired distribution.
\end{proof}

Finally we relate a colouring of a product to a colouring of one of its factors. For a fixed coordinate \(i\), colour each vertex \(v\) of that factor by a most frequent colour on the fibre \(x_i=v\), with respect to the product measure. The error of this rule is the measure of the product vertices whose colour it predicts incorrectly.
On a product probability space write \(f_T\) for the orthogonal component of \(f\) centred in exactly the coordinates \(T\), and \(\tail(f)=\sum_{|T|\ge2}\|f_T\|_2^2\).
\begin{lemma}[A single coordinate from small tails]\label{lem:one-coordinate}
If colour indicators \(f_c\) partition a product probability space, have masses at most \(a\), and \(\sum_c\tail(f_c)=E\), then one coordinate's modal classifier has error \(O(a+E)\).
\end{lemma}
\begin{proof}
Write \(F\) for the vector of indicators. For a bipartition of the coordinates, let \(U,V\) be the centred conditional expectations of \(F\) on the two blocks, and put \(S=\E F+U+V\). Since \(FF^{\mathsf T}=\operatorname{diag}F\), \(\|S\|\le3\), and \(\E\|F-S\|^2\le E\), we have
\(\E\|SS^{\mathsf T}-\operatorname{diag}S\|_{\mathrm F}^{2}=O(E)\),
where \(\|\cdot\|_{\mathrm F}\) is the Frobenius norm. Its component centred in both blocks is \(UV^{\mathsf T}+VU^{\mathsf T}\), of squared norm at least \(2\|U\|^2\|V\|^2\). Independence gives \(\E\|U\|^2\,\E\|V\|^2=O(E)\). Thus the product of the first-level weights in any two complementary blocks is \(O(E)\). The total is \(1-O(a+E)\), so one coordinate carries all but \(O(a+E)\) of it. Predicting the most frequent colour on each fibre makes an error bounded by the summed conditional variance, since \(1-\max p_c\le1-\sum_c p_c^2\).
\end{proof}

We shall use the \((\ell,0)\)-tree suspension of Ai, Ding, Liu and
Yang~\cite[Definitions~2.1 and~2.3]{ADLY} with a finite forbidden
family \(\mathcal A\) as input. At each root and at each node of an
attached choice tree, choose a member of \(\mathcal A\) separately.
Let \(\mathcal S_\ell(\mathcal A)\) consist of all the resulting \(k\)-graphs.
This family is finite since the depth is fixed and every member of
\(\mathcal A\) has at most \(M=\max_{A\in\mathcal A}|V(A)|\) vertices.
The density formula for a single input hypergraph was proved
in~\cite[Theorem~2.6]{ADLY}.

\begin{lemma}[Tree suspensions of finite families]\label{lem:family-suspension}
Let \(1\le\ell<k\) and let \(\mathcal A\) be a nonempty finite family of
\(k\)-graphs, each with at least one edge. Put \(\beta=\pi_\ell(\mathcal A)<1\),
\(s=k-\ell\) and \(c=(1-\beta)^{-1/s}\). Then
\(\pi_\ell(\mathcal S_\ell(\mathcal A))=1-(1+c)^{-s}\).
Deleting any independent set from a member of \(\mathcal S_\ell(\mathcal A)\)
leaves a copy of a member of \(\mathcal A\).
\end{lemma}
\begin{proof}
Let \(I\) be independent in a suspension. If \(I\) meets every root,
choose one of its vertices in each root and follow the tree attached
to this transversal. At each node met by \(I\), choose a vertex in
\(I\) and continue along its branch. A full branch would give a
connecting \(k\)-edge in \(I\). Thus some root or node is disjoint
from \(I\), proving the deletion assertion.

Put \(\theta=1/(1+c)\) and \(\gamma=1-\theta\).
Join an independent class of proportion \(\theta\) to a near-extremal
\(\mathcal A\)-free core of proportion \(\gamma\), including all crossing
edges. The deletion assertion excludes every suspension. An \(\ell\)-set
in the independent class loses a proportion \(\theta^s+o(1)\) of its
extensions, and one in the core loses at most
\((1-\beta)\gamma^s+o(1)=\theta^s+o(1)\). A mixed set has every extension.

For the reverse inequality, take a host whose minimum \(\ell\)-degree
ratio exceeds \(1-\theta^s\) by a fixed margin. On any set of at least
\(\gamma n-O(1)\) vertices the induced ratio exceeds \(\beta\) by a fixed
margin, since \((1-\beta)\gamma^s=\theta^s\). Such a set therefore contains
some member of \(\mathcal A\). First embed the \(\ell\) disjoint roots.
In the link of each root transversal, use the degree count in
\cite[Lemmas~2.7 and~2.8]{ADLY}. In a \(j\)-uniform link of density
above \(1-\theta^j\) by a fixed margin, for \(j\ge2\) at least
\(\gamma n-O(1)\) vertices have links of density above
\(1-\theta^{j-1}\) by a fixed margin. Embed a member in this vertex
set and recurse at each of its vertices. For \(j=1\) the link itself
is a set of size at least \(\gamma n-O(1)\). The depth and \(M\) bound
the total number of used vertices, so deleting them preserves the margins.
The members chosen in this embedding form one of the suspensions in
\(\mathcal S_\ell(\mathcal A)\).
\end{proof}

\section{Ordinary densities}\label{sec:ordinary-cases}

We prove the ordinary-density case by constructing an \(r\)-graph whose extremal blow-ups have edge density \(\beta\) with \(\Q((1-\beta)^{-1/(r-1)})=K\). We first build the pattern for triples and then for larger uniformities. We next adjoin independent vertex classes to its blow-ups and construct a single forbidden hypergraph with the resulting density. Increasing the number of added classes makes these densities tend to one.

\subsection{The extremal patterns}\label{sec:ordinary-classes}
Let \(J\) be a finite \(r\)-graph and let \(x_v\) be the proportion of vertices in the blow-up class corresponding to \(v\in V(J)\). The limiting edge density of this blow-up is \(P_J(\mathbf x)=r!\sum_{e\in E(J)}\prod_{v\in e}x_v\). Write \(\lambda(P_J)\) for its maximum over \(x_v\ge0\) with \(\sum_v x_v=1\). We shall choose \(J\) so that this maximum gives the desired field parameter.

The starting point for the choice of class proportions is the following fact about cubics.
\begin{theorem}\label{ch:maxima}
Every real number field is generated by the maximum of a rational homogeneous cubic on a simplex. The maximum can be unique and interior and the multinomially normalized coefficients can lie in \((1/3,2/3)\).
\end{theorem}
\begin{proof}
Index the coordinates \(t_j\) by all embeddings \(K\hookrightarrow\mathbb C\),
starting with the given real embedding. On \(K_\infty\), conjugate
embeddings have conjugate coordinates, so the trace hyperplane is
\(\sum_j t_j=1\). The idempotent of the first embedding is
\(e=(1,0,\ldots,0)\). For a primitive element \(\alpha\) with minimal
polynomial \(m\), its coordinates in the power basis are the coefficients
of \(m(T)/((T-\alpha)m'(\alpha))\), which lie in \(K\). Given \(b\in K\), the rational polynomial \(P_b=\operatorname{Tr}(b(3x^2-2x^3))\) has value \(b_1\) and zero gradient at \(e\). Choose \(b\) primitive with \(b_1>b_j\) for every other real embedding. On the trace hyperplane its quadratic part at \(e\) is
\[
3\sum_{j>1}(b_j-b_1)t_j^2-6b_1\sum_{1<i<j}t_it_j .
\]
Add a homogeneous cubic vanishing to first order at the conjugate idempotents. This space is rational (expand the jets at \(e\) in a basis of \(K\)); over \(\mathbb C\) it consists of distinct-index monomials. Choose the \(t_1t_it_j\) coefficients near \(6b_1\), except for conjugate pairs where they are near \(6b_1-M\), \(M>6\max_j|b_j-b_1|\). The real Hessian is negative definite, since each nonreal pair contributes \(6\Re((b_j-b_1)t_j^2)-M|t_j|^2\); rational approximation in the jet space preserves this. A small rational simplex with \(e\) interior then gives a strictly concave objective. Barycentric homogenization yields rational coefficients, which approach \(b_1\) as the simplex shrinks. A rational affine rescaling gives the requested range. The same argument works for any prescribed \(b\) satisfying the real ordering, and any finite list of other specified stationary values is obtained on the same simplex by their uncorrected \(P_b\)'s.
\end{proof}

We first realize such a cubic by triples of vertices.
\begin{theorem}\label{ch:ordinarythree}
For every number field \(K\subset\R\) there are single forbidden \(3\)-graphs whose ordinary Tur\'an densities generate \(K\) and approach one.
\end{theorem}
For the construction it suffices to find a finite pattern of edge density \(\pi_*\) for which \((1-\pi_*)^{-1/2}\) generates \(K\). We pass from this pattern to a single forbidden \(3\)-graph in Section~\ref{sec:ordinary-single}.

The next estimate will let us prescribe the numbers of edges of each
class profile. The shift adds the same element of \(\mathbb Z/N\mathbb Z\)
to every index in a tuple. In each profile we use full shift orbits
representing edges with distinct vertices. The tensor is centred by
subtracting its expectation under independent selection of these orbits.
\begin{lemma}[Cyclic orbit estimate]\label{lem:cyclic-orbit}
For fixed \(r\), one can choose a prescribed number of eligible full shift orbits on \(\Z/N\Z\) with centred ordered tensor of injective Euclidean norm \(O_r(N^{3/4}\log N)\).
\end{lemma}
\begin{proof}
Select each eligible orbit independently with probability chosen to give the required expected count. For unit test vectors and orbit sums \(b_O\), unitary Fourier orthogonality gives
\[
 \sum_O |b_O|^2
 =N\sum_{\xi_1+\cdots+\xi_r=0}\prod_j|\widehat v_j(\xi_j)|^2
 \le N\min_j\|\widehat v_j\|_\infty^2
\]
over ordered Cartesian orbits. Restrict to eligible orbits and group at most \(r!\) ordered orbits per unordered one, costing only a constant by Cauchy--Schwarz. Split the Fourier vectors at \(N^{-1/4}\). A term with a small part has squared ranges \(O_r(\sqrt N)\) and an \(N^{-r}\)-net of logarithmic size \(O_r(N\log N)\).

The all-large term has squared ranges \(O_r(N)\) but supports of size \(\sqrt N\), hence logarithmic net size \(O_r(\sqrt N\log N)\). Hoeffding and the union bound give the stated norm with failure \(\exp(-\Omega_r(\sqrt N\log^2N))\); the Frobenius bound \(O_r(N^{r/2})\) transfers from the nets. The desired count is a binomial mode, of probability at least \(1/(|\mathcal O|+1)\), where \(\mathcal O\) is the eligible orbit set.
\end{proof}

\begin{proof}[Proof of Theorem~\ref{ch:ordinarythree}]
Choose \(w\) as in Lemma~\ref{lem:choose-generator}, with selected value positive and near \(1/6\) and the other real values larger and positive. Our choice of triples is guided by one requirement:
after the vertex weights are optimized, \(\sqrt{1-\pi_*}\) should be a
nonconstant rational fractional-linear function of \(w\). We first encode
\(w\) in the triples on the data vertices. The triples involving the roots
will supply the other coefficients of the cubic in the total data weight.

By the construction in the proof of Theorem~\ref{ch:maxima}, on a common small simplex there are rational cubics \(G_0,G_1\) stationary at the same interior point, with \(G_0\) strictly concave and respective stationary values \(1-18w^2\), \(1+12w\). For a sufficiently large divisible square \(L\), \(G_L=G_0-G_1/L\) is still strictly concave and has maximum
\begin{equation}
\beta=1-18w^2-(1+12w)/L. \label{eq:ordinarythree-2}
\end{equation}
Its multinomially normalized coefficients lie in \((.4,.6)\) and have denominators dividing \(D_0L\), with \(D_0\) fixed.

Replace each of the \(d\) simplex coordinates by a class of \(L\) vertices indexed by \(\Z/L\Z\). We choose triples in whole orbits under simultaneous translation of the three indices. When the weights within each class are equal, an orbit of full size contributes \(1/L^2,2/L^2,6/L^2\) to the normalized coefficient according as the class multiplicities of the triple are \(111,21,3\). We use only triples of distinct vertices and discard the short orbits of equally spaced triples within a class. Taking \(6D_0\mid L\) makes the number of orbits required for each coefficient an integer, and the proportion required stays in a compact subinterval of \((0,1)\).

Choose the orbits in each profile by Lemma~\ref{lem:cyclic-orbit}. Let \(P_D\) be the edge polynomial of the triples just chosen. Give these vertices weights \(z\) of total mass \(t\). Write \(\bar z\) for the weights obtained by distributing the mass of each class equally among its \(L\) vertices, and put \(\delta=z-\bar z\). At \(\bar z\) the centred tensor has value zero, and its first variation within each class vanishes by translation invariance. The expected polynomial differs from the constant-profile polynomial only on repeated vertices and short orbits within a class. In directions summing to zero within each class, these terms have Hessian norm \(O_d(t)\), since the short orbits have bounded vertex degree. Thus, with \(H=O_d(L^{3/4}\log L)\),
\begin{equation}
 |P_D(z)-P_D(\bar z)|\ll_d(t+Ht/\sqrt L+H\|\delta\|)\|\delta\|^2. \label{eq:ordinarythree-5}
\end{equation}
At total mass one, we also have \(\lambda(P_D)\le1-\eta\) for \(\eta\asymp H^{-2/3}\). Indeed, a larger value would give \(\|z\|^2<\eta\), since a repeated vertex cannot lie in a triple of distinct vertices. The expected value is at most \(.75\), and the centred contribution is at most \(H\eta^{3/2}<.1\) if the constant in \(\eta\) is sufficiently small.

Adjoin \(q=3(L-2)\) vertices, which we call roots. On these vertices take all triples except ten disjoint perfect matchings. To obtain the matchings, use \((i,i+a,i+2a)\) on three copies of \(\Z/(L-2)\) for ten choices of \(a\). Fix a cycle through the roots. A triple with two roots and one data vertex is an edge unless the roots are consecutive on this cycle. For triples containing just one root \(r\), permute the fine indices in each data part by \(\sigma_{i,r}\), and exclude data pairs whose permuted indices agree. 

Choose permutations with \(q^{-1}\sum_r\|\sum_iP_{\sigma_{i,r}}\delta_i\|^2\ge\|\delta\|^2/2\) on centred data vectors: the independent random forms average to the identity, have norm at most \(d\), and matrix Chernoff~\cite{Tropp} applies since \(q\asymp L\). For the moment distribute the root mass uniformly and distribute the mass of
each data class uniformly among its fine vertices. Writing \(s=1-t\), and
optimizing the data class proportions at fixed \(t\), we obtain
\begin{equation}
 F_\beta(t)=a s^3+3b s^2t+3cst^2+\beta t^3,\quad
 (a,b,c)=(1-3/q-18/q^2,1-3/q,1-1/L). \label{eq:ordinarythree-8}
\end{equation}
Among three root vertices the complete \(3\)-graph contributes \(1-3/q+2/q^2\), and removing the ten matchings costs \(20/q^2\). With two root vertices, repetitions and the edges of the root cycle account for the loss \(3/q\). With one root vertex, each common index class of data vertices has weight \(t/L\), giving the loss \(1/L\).

We next show that every maximizing weighting has these two uniformity
properties. Distributing all the weight uniformly among the roots attains \(a\). Thus any maximizing weighting satisfies \(\eta t^3\le1-a\), so \(t=o(1)\). Counting repetitions among the roots gives \(\max_r y_r=O(q^{-1/2})\) for the root weights \(y\). Put \(e=y-s{\bf1}/q\). The root-only quadratic loss \(3s\|e\|^2\) dominates the \(O(\max y_r)\|e\|^2\) Hessian contributions of \(\sum y_r^3\) and the bounded-degree deletion; first variations vanish. The two-root term costs only \(O(t)\|e\|^2\) above its uniform value. The one-root term is
\(3cst^2-3\sum_r y_r\|\sum_iP_{\sigma_{i,r}}\delta_i\|^2\).
With \(u=t/(1-t)\) and \(v=1-\eta\), the numerator of \(F_v-a\) over \((1+u)^3\) is \(3(b-a)u+3(c-a)u^2+(v-a)u^3\). Since \(b-a,c-a=O(q^{-2})\), \(v-a\le-\eta/2\), this gives \(t=O(1/(q\sqrt\eta))\) and \(F_v-a=O(t/q^2)\). Counting the loss on triples within the root class now yields \(\|e\|^2=O(t/q^2)\), hence \(y_r\ge s/(2q)\). The permutation inequality likewise gives \(\|\delta\|^2=O(t/q^2)\). Thus the coefficient in \eqref{eq:ordinarythree-5} is \(o(1)\), smaller than the one-root quadratic loss. Every maximizer is therefore uniform in the fine parts and at the roots.

It remains to choose the total data weight. In terms of \(u=t/(1-t)\), the positive critical point solves \((b-a)+2(c-b)u+(\beta-c)u^2=0\), and is unique because \(b>a,\ \beta<c\). Equations \eqref{eq:ordinarythree-2} and \eqref{eq:ordinarythree-8} give \(u=1/(qw)\). The optimized value is \(\pi_*=(a+2bu+cu^2)/(1+u)^2\), and
\[
\sqrt{1-\pi_*}=\frac{3\sqrt L/q+u/\sqrt L}{1+u},
\]
a nonconstant rational fractional-linear function of \(u\) since \(q\ne3L\) and \(L\) is square.

Blow up the pattern in the maximizing proportions to obtain \(3\)-graphs of limiting density \(\pi_*\). Every maximizer gives positive weight to every vertex, so the pattern is minimal. Pikhurko's theorem~\cite{P} now gives a finite forbidden family with this same density. The extremal blow-ups have bounded chromatic number, and \(c=(1-\pi_*)^{-1/2}\) generates \(K\).
For this finite pattern, the conversion in Theorem~\ref{thm:ordinary-conversion}
gives singleton densities \(1-(q+c)^{-2}\). They approach one and
generate \(K\) for all sufficiently large \(q\) by
Lemma~\ref{lem:field-transfers}.
\end{proof}

For larger uniformities we use the same distinction between root and data vertices, choosing the edges according to how many vertices they have in each class.
\begin{theorem}\label{thm:ordinary-fields}
For every \(r\ge3\) and every number field \(K\subset\R\), single forbidden \(r\)-graphs have ordinary Tur\'an densities generating \(K\) arbitrarily close to one.
\end{theorem}

\begin{proof}
The case \(r=3\) is Theorem~\ref{ch:ordinarythree}.
For \(r\ge4\) put \(s=r-1\), fix rational \(\theta\in(0,1)\) with \(\theta^s>1/2\), and choose a primitive \(w\) positive and near \((2(1-\theta))^{-1/s}\), with larger positive real conjugates. Apply the construction in the proof of Theorem~\ref{ch:maxima} to \(B=1-(1-\theta)w^s\), then multiply by \((\sum x_i)^{r-3}\), obtaining a rational degree-\(r\) polynomial \(G_0\) with unique interior maximum \(B\) and normalized coefficients in \((.4,.6)\) (homogenization averages the cubic coefficients).

Take a single root class and, for each variable of \(G_0\), a separate data
class. Give each class \(N\) vertices. The profile of an \(r\)-set records its number of vertices in each class. Take a divisible integer \(A\) with \(A^s\sim2\binom r2 N^{s-1}\) and put \(a=A/N,\ b=\theta a\). Prescribe coefficients \(1-a^{s-j}b^j\) on profiles with \(j\le s\) data vertices, and use \(G_N=G_0-\theta b^s(\sum x_i)^r\) on all-data profiles, of maximum \(\beta=B-\theta b^s\). A full shift orbit of multiplicities \(m_i\) contributes \(\prod m_i!/N^{r-1}\); divisibility by \(r!\) and the coefficient denominators gives the exact counts. Repetitions cost at most \(\binom r2/N\), and short orbits \(O_r(N^{-r/2})\) (in a shift stabilizer every position belongs to a run of length at least two), so \(\theta^s>1/2\) ensures sufficient eligible orbits.

Choose the translation orbits for each edge profile by Lemma~\ref{lem:cyclic-orbit}, so that its centred tensor has norm \(H=O_r(N^{3/4}\log N)\). At a maximizer \(\|x\|^2=O(1/N)\) by the missing repeated-label probability and the uniform root-only value. The all-data cap is \(1-N^{-1/2}\): otherwise at mass one simplicity gives norm at most \(N^{-1/4}\), so the centred error \(O(HN^{-r/4})=o(1)\) cannot close the constant gap in expectation. Thus the total weight on the data vertices tends to zero. Now average the weights within each class and write \(\delta\) for the difference from these averages. Translation invariance makes all first-order terms in \(\delta\) vanish. The terms with at least two deviations give an error
\[
 O((H+1)N^{-(r-2)/2}\|\delta\|_2^2)=o(\|\delta\|_2^2)
\]
by the injective norm bound. Here the short-orbit norm is \(O_r(1)\): split into shift cycles, each contracted sum of products of at least two Euclidean unit vectors being bounded by H\"older, and sum over the finitely many cycle types. In the remaining complete-simple expectations elementary symmetric polynomials increase under averaging. On the root-only profile the gain is \(\Omega(\|\delta_{\rm root}\|^2)\): pair averaging gains \((x-y)^2e_{r-2}(\text{others})/4\), with \(e_{r-2}\) bounded below by root mass \(1-o(1)\) and maximum weight \(o(1)\). For each data part use its profile with exactly two data vertices, whose \(e_2\) gain is half its variance loss, with the same root factor. Consequently every maximizing weighting is constant on each class.

We may therefore optimize first within the data classes and then over the total root mass \(x\). Writing \(t=1-x\), the resulting polynomial is \(P=\sum_{j=0}^s\binom rj(1-a^{s-j}b^j)x^{r-j}t^j+\beta t^r\).
Putting \(L=ax+bt\) gives \(P_x/r=1-L^s,\ P_t/r=1-\theta L^s-(1-\theta)(wt)^s\) on the simplex.

Hence the unique maximum is at \(t=a/(w+a-b)\); Euler's identity gives \(1-a_*^s\) with \(a_*=aw/(w+a-b)\), an invertible rational fractional-linear transform of \(w\).

Taking blow-ups in these proportions gives the matching lower bound. Each vertex of the pattern has positive weight at every maximum, so Pikhurko's theorem again gives a finite forbidden family with the same ordinary density.
Here \(c=1/a_*\) generates \(K\). Applying
Theorem~\ref{thm:ordinary-conversion} to the pattern gives the singleton
density \(1-(q+c)^{-s}\). This tends to one and, by
Lemma~\ref{lem:field-transfers} with \(R(X)=1-X^{-s}\), generates
\(K\) for all sufficiently large \(q\).
\end{proof}

\subsection{A single forbidden hypergraph}\label{sec:ordinary-single}\label{ch:product}
Let \(\D\) consist of the subhypergraphs of blow-ups of the pattern just
constructed. Write \(\beta\) for its density and \(c=(1-\beta)^{-1/(r-1)}\),
so that \(\Q(c)=K\). Adjoin \(q\) independent vertex classes to a core in
\(\D\), allowing every edge meeting two or more components. We call the
added classes \emph{pads}. Balancing the missing edges within the
components gives proportion \(1/(q+c)\) to each pad and \(c/(q+c)\)
to the core, and hence density \(\alpha=1-(q+c)^{-(r-1)}\).

To realize this density with one forbidden hypergraph we use a product.
In \(B_1\times B_2\), the vertices \(\{v\}\times V(B_2)\) form the
\emph{row} indexed by \(v\in V(B_1)\), and \(w\in V(B_2)\) specifies a
position in that row. In a product with more factors, fix one factor
\(B_j\) to index the rows. An \emph{outer position} is a choice of one
vertex in each of the other factors. We call a list of positions
\emph{complete} if its entries in each outer factor are distinct and span
a clique. An edge of \(B_j\) then lifts to an edge of the product by
placing its vertices at distinct positions on the list. A \emph{coarse
colour} records a pad or the core. Within the core the representation
may give a vertex a further label.

\begin{theorem}\label{thm:ordinary-conversion}
Let \(r=s+1\ge3\), and let \(\D\) be a represented class of weakly
\(p\)-colourable \(r\)-graphs of density \(\beta<1\). Assume that its
representations restrict to subhypergraphs and that it has a finite
exact family of witnesses outside \(\D\).
Suppose that, given \(L\ge Cb\) complete outer positions, a representation
on the rows using at most \(b\) coarse colours per row induces a
representation of the factor indexing the rows in the same padded
class. Then \(1-(q+(1-\beta)^{-1/s})^{-s}\) is a singleton density for
all sufficiently large \(q\), depending initially only on \(s,p,C\).
\end{theorem}
For a finite pattern, give each row its most frequent coarse colour. In
a core row keep its most frequent fine label among the positions of
that coarse colour. For \(L/b\) sufficiently large in terms of the
fixed pattern and \(r\), each chosen label occurs in at least \(r\)
positions. The vertices of any edge can therefore be placed at distinct
positions carrying their chosen labels.

Fix \(q\), and put \(\alpha=1-(q+(1-\beta)^{-1/s})^{-s}\).
Starting with the exact witness family for \(\D\), apply
Lemma~\ref{lem:family-suspension} with \(\ell=1\) a total of \(q\) times.
Each step increases \((1-\beta')^{-1/s}\) by one, where \(\beta'\) is
the current density. The resulting finite family \(\mathcal B_q\)
has density \(\alpha\). Its members lie outside the padded class, since
deleting the \(q\) independent classes in turn would leave an original
witness in the core. Put \(\zeta=(1-\alpha)^{1+1/s}/10\).
For a clique order \(L\) to be chosen below, take \(m\) sufficiently large
and divisible by \(L\). Let \(P\) be the product of all \(m\)-vertex
\(r\)-graphs that contain a member of \(\mathcal B_q\) and have minimum
vertex-degree ratio at least \(\alpha-\zeta\).

The coordinate projections will give the upper bound. For the reverse
inequality we shall use a colouring of the product to colour one of
its factors. We first estimate correlations among the positions of a
product edge.
\begin{lemma}[Product correlations]\label{lem:product-correlations}
In each factor use an \(s\)-tuple law of density \(M_i\) relative to \(\sigma_i^s\), with all marginals \(\sigma_i\) and \(\|M_i-1\|_2\le d\). Sample independently across factors. For indicators in the sampled positions the correlation differs from the product of means by \(O_s(d)\), uniformly in the number of factors. There is a constant \(d_0=d_0(s)>0\) such that, when \(0\le d\le d_0\), the error for a common indicator of mean \(a\) and first-level weight \(V\) is at most
\begin{equation}
dV^{s/2}+O_s(da^{(s+1)/2}+d^2a^{s/2}).\label{eq:product-1}
\end{equation}
\end{lemma}
\begin{proof}
Decompose orthogonally by centred coordinates. Each active factor uses a set of at least two positions, contributes norm at most \(d\), and forces those positions' indicator components to be centred there. For one active factor, sum over its index by Cauchy--Schwarz on two positions. For a common indicator and \(j\) active positions the sharper bound is \(da^{s-j}V^{j/2}\). For at least two active factors majorize by a tensor power of the nonnegative bit tensor equal to \(1\) at \(0^s\), zero at weight one, and \(d_0\) elsewhere. Its injective norm is \(1\) for small fixed \(d_0(s)\): on inputs \((x_j,y_j)\ge0\) of unit norm the constant term is at most \(1-\sum y_j^2/(2s)\), and the remainder is \(O_s(d_0\sum y_j^2)\). This bound tensorizes by applying it to the slice norms. For \(0\le d\le d_0\), every term with at least two active factors is at most \((d/d_0)^2\) times its majorant at \(d_0\). Their sum is therefore \(O_s(d^2)\prod_j\|f_j\|_2\). For \(d>d_0\) the coarse estimate follows because the correlation error of indicators is at most \(1\le d/d_0\).
\end{proof}

The following consequence will also be used for the links of vertices
in the four-uniform construction.
\begin{lemma}[Independent sets in products]\label{lem:ordinary-product-tail}
Fix \(s\ge2\). Let the factors of a finite categorical product be \((s+1)\)-graphs
of order \(m\), each with minimum vertex-degree ratio at least \(1-\mu\).
Put \(\lambda=\mu^{1/s}\), where \(\mu>0\) is sufficiently small and
\(m\) is sufficiently large in terms of \(s\) and \(\mu\).
If \(f\) is the indicator of an independent set of density \(a\) in
the product, with the uniform measure on its vertices, then
\[
 a\le\lambda+O_s(\lambda^{3/2}),\qquad
 \tail(f)=O_s\bigl((\lambda-a)_++\lambda^{3/2}\bigr).
\]
The bounds are uniform in the number of factors.
\end{lemma}

\begin{proof}
For the empty independent set both bounds are immediate. We apply Lemma~\ref{lem:product-correlations} to links of vertices in the factors. For \(m\)-vertex factors of minimum vertex-degree ratio at least \(1-\mu\), put \(\lambda=\mu^{1/s},\ \tau=\sqrt\lambda,\ \nu=\mu/\tau+1/m\), with \(m\) large. Call a pair good if its degree ratio is at least \(1-\tau\); averaging the missing edges through any vertex shows that all but an \(O_s(\nu)\) fraction of the other vertices make a good pair with it.

We use the following normalization. On a finite probability space of law \(\sigma\), let \(D\) be a symmetric missing-\(s\)-tuple indicator, \(r(x)=\E[D\mid X_1=x]\), \(\eta=\|r\|_\infty\), and let \(K\) be the kernel \(\E[D\mid X_1=x,X_2=y]\). Solve
\[
 (I+(s-1)J-\operatorname{diag}(r)-(s-1)K)u=r,\qquad Ju=\E_\sigma u .
\]
The inverse of \(I+(s-1)J\) is bounded, and the perturbation has norm \(O_s(\eta)\) in both sup norm and \(L^2\), by the row bounds and symmetry. Thus \(\|u\|_\infty=O_s(\eta)\), \(\|u\|_2^2=O_s(\eta\,\E D)\). For small \(\eta\), \(M=(1-D)(1+\sum_j u(x_j))\) is a nonnegative density with marginals \(\sigma\), by the displayed row equation. On missing tuples its squared deviation from one is one. On the remaining tuples the mean squared deviation is at most \(s^2\|u\|_2^2\). Consequently \(\|M-1\|_2^2\le(1+O_s(\eta))\E D\).

Apply the normalization first to the good-pair graph, with uniform measure. Its stochastic kernel \(B_i\) has entries at most \((1+O_s(\nu))/m\) and centred norm \(O_s(\nu)\), by subtracting it from a constant kernel of that height. For a vertex \(z\) of the \(i\)th factor, take \(\sigma=B_i(z,\cdot)\) and apply the same normalization to the \(s\)-edges of its link. The row defect is \(O_s(\tau+\nu)\) by goodness, and total defect at most \(\mu(1+O_s(\nu))+O_s(1/m)\). We thus obtain \(d=\|M-1\|_2\le\sqrt\mu(1+O_s(\sqrt\lambda))\). Taking \(\mu\) sufficiently small and then \(m\) sufficiently large puts \(d\) in the range of \eqref{eq:product-1}.

Let \(f\) be the indicator of an independent set in the product and let \(z\) belong to that set. No edge through \(z\) has all its other vertices in the set, so the corresponding correlation in the product of links is zero. With conditional mean \(a_z\) and first-level weight \(V_z\), \eqref{eq:product-1} gives \(a_z^{s/2}\le d(V_z/a_z)^{s/2}+O_s(d\sqrt{a_z}+d^2)\) if \(a_z>0\). First \(a_z=O_s(\lambda)\), then \(a_z\le\lambda+O_s(\lambda^{3/2})\) and conditional tail \(O_s((\lambda-a_z)_++\lambda^{3/2})\), by inversion for \(a_z\ge\lambda/2\) and trivially otherwise.

We transfer the estimates from the link of a particular root to the
uniform measure on the product. The product kernel has centred norm \(O_s(\nu)\), so \(\E(|a_z-a|\mid f(z)=1)=O_s(\nu)\). Two independent partners \(X,Y\) of a uniform root have joint density in each factor at least \(1-O_s(\nu^2)\), from averaging \((H-b_z(x))(H-b_z(y))\ge0\) for \(b_z=mB_i(z,\cdot)\le H=1+O_s(\nu)\). In one factor, the density of the pair \((X,Y)\) conditional on the root,
relative to its unconditional law, is at most
\(H^2/(1-O_s(\nu^2))=1+O_s(\nu)\). After subtracting one, the kernel
has mean zero in each row and column. Its positive part is bounded by
\(O_s(\nu)\), so its absolute row and column integrals are \(O_s(\nu)\).
Schur's test gives norm \(O_s(\nu)\) on mean-zero functions. The product
channel has the same bound, since each nonconstant tensor component has
a mean-zero factor and the other factors are contractions.

To estimate the tail under uniform measure we use a statistic whose
expectation vanishes on the constant and single-coordinate components.
Apply the channel bound to
\(f(X)^2+3f(X)f(Y)-\sum_{j=1}^4 f(X)f(X_{I_j}Y_{\overline I_j})\),
where each coordinate is assigned independently and uniformly to one of
the four blocks \(I_j\). The statistic has norm \(O(\sqrt a)\), so
conditioning on \(f(z)=1\) changes its expectation by \(O_s(\nu)\).
Its conditional expectation at a root is at most the conditional tail.

For the unconditional expectation, diagonalize each \(B_i^2\) on
mean-zero functions and adjoin the constant eigenfunction. On a product
basis vector with nonconstant coordinates \(J\), let \(\rho_i\) be the
corresponding eigenvalue of \(B_i^2\) for \(i\in J\). The multiplier is
\[
  1+3\prod_{i\in J}\rho_i
    -4\prod_{i\in J}\frac{1+3\rho_i}{4}.
\]
It vanishes for \(|J|\le1\). For \(J=\{i,j\}\) it equals
\(\tfrac34(1-\rho_i)(1-\rho_j)\). Since \(0\le\rho_i=O_s(\nu^2)\),
the multiplier is at least \(1/2\) on every basis vector with
\(|J|\ge2\) for small \(\nu\). Summing over the product basis bounds
the uniform tail by twice the unconditional expectation. Since \(\nu=O_s(\lambda^{3/2})\), the uniform bounds are \(a\le\lambda+O_s(\lambda^{3/2})\), \(\tail(f)=O_s((\lambda-a)_++\lambda^{3/2})\).
\end{proof}

\begin{proof}[Proof of Theorem~\ref{thm:ordinary-conversion}]
The upper bound now follows from Lemma~\ref{lem:product-upper}. In a
host whose minimum vertex-degree ratio exceeds \(\alpha\) by a fixed
margin, keep a bounded witness and sample the remaining vertices of an
\(m\)-set. Concentration gives an induced subhypergraph with the
required minimum degree, so the host contains a factor. The coordinate
projections and \(\pi_1=\pi\) give \(\pi(P)\le\alpha\).

For the lower bound suppose \(P\) has a representation in the padded
class. A product with just one factor already contains its forbidden
witness, so assume there are at least two factors.
For the product \(P\) constructed above take \(\mu=1-\alpha+\zeta\).
The representation induces a weak colouring of \(P\) with \(t=p+q\)
colours, where \(t\mu^{1/s}-1=O(\mu^{1/s})\).
Write \(\lambda=\mu^{1/s}\) and \(\nu=\mu/\sqrt{\lambda}+1/m\).
We use the neighbour kernels and link distributions from the proof of
Lemma~\ref{lem:ordinary-product-tail}.
Apply Lemma~\ref{lem:ordinary-product-tail} to the \(t\) independent colour classes. When \(t\lambda-1=O_p(\lambda)\), the sum of their tails is \(O_{s,p}(\sqrt\lambda)\), since the masses sum to one. Lemma~\ref{lem:one-coordinate} gives a coordinate whose most frequent colour on each fibre has this error.

Choose \(\eta>0\) with \(2(1-2\eta)^s>1\), and index the rows by the vertices of this factor. Apart from \(O_{s,p}(\lambda^{-1/2})\) exceptional colours, each colour has a set of \((1-\eta)\lambda m\) rows on which it occurs with probability at least \(1-1/(10s)\) over the outer positions. These sets of rows are disjoint. The exceptional colours are accounted for by the mass deficits or the error of the one-coordinate colouring. Average over the outer neighbour kernels and apply Chebyshev's and Markov's inequalities. Outside a set of outer positions of measure \(O_{s,p}(t\nu^2)\), each such colour has at least \((1-2\eta)\lambda m\) rows on which its averaged probability exceeds \(1-1/(2s)\). Call the remaining outer positions retained.

Suppose a row \(x\) has this colour at a retained position. No \(s\)-set of the corresponding rows can form an edge with \(x\). Otherwise the link distributions in the outer factors and a union bound would yield a monochromatic edge in the product. Two nonexceptional colours on \(x\) would therefore give more missing edges through \(x\) than its degree permits, by the choice of \(\eta\). Hence \(x\) uses only \(O_{s,p}(\lambda^{-1/2})\) colours over all retained positions.

Take \(L\asymp\lambda^{-1/2}\), large enough to decode. By Lemma~\ref{lem:uniform-cliques}, for \(L\mid m\) each outer factor has a uniform-marginal clique law, since \(\mu L^s=O(\lambda^{s/2})\). Independent aligned choices across factors let all \(L\) positions be retained since \(Lt\nu^2=o(1)\), so the row-colouring rule gives a representation of the chosen factor.

A representation of \(P\) in the padded class would therefore induce a
representation of one factor, contradicting the witness it contains.
The padded class avoids \(P\) and gives \(\pi(P)\ge\alpha\), as required.
\end{proof}

\section{Codegree densities}\label{ch:cg}\label{app:cgproof}

For codegree we first construct a finite forbidden family of \(3\)-graphs by prescribing the neighbourhoods of pairs in a cyclic partition. Its extremal density can be chosen to generate \(K\), with control over all its conjugates. We then lift the construction one uniformity at a time by prescribing vertex links. Finally we join the extremal examples to independent classes and use a product to obtain a single forbidden hypergraph with codegree density tending to one.

\begin{theorem}\label{cg:main}
Fix \(k\ge3\) and a number field \(K\subset\R\). For any \(u\in(0,1)\) and \(\varepsilon>0\), there is a finite-family codegree density \(\gamma\) generating \(K\), with \(|\sigma(\gamma)-u|<\varepsilon\) for every embedding \(\sigma:K\to\mathbb C\). There are also singleton codegree densities generating \(K\) arbitrarily close to one. In both statements the density can be chosen to be the least of its real conjugates.
\end{theorem}

\subsection{The extremal class}\label{sec:codegree-classes}
We first construct a coloured class of \(3\)-graphs. Number its independent
vertex classes \(V_j\) cyclically modulo \(p\). For a reduced rational
\(A_0/c\in(0,1)\), choose \(p,b\) with \(A_0p-bc=1\) and \(p\) large, so
\(b/p\) is close to the target \(A_0/c\). If \(u,v\in V_j\) are distinct,
their neighbours may only lie in the classes indexed by
\(W_j=\{j+1,\ldots,j+b\}\pmod p\). Thus, writing \(x_a=|V_a|/n\) and
taking \({\bf1}_{W_j}\) as row \(j\) of a matrix \(C\), we have
\[
 \frac{\delta_{\co}(H)}{n}\le\frac{|N_H(u,v)|}{n}
 \le\sum_{a\in W_j}x_a=(Cx)_j.
\]
These are the codegree bounds from pairs within a class. For selected
pairs in different classes we shall restrict the graph links of
individual vertices. The resulting losses will determine how far the
density falls below \(b/p\).

We first record how the window rule produces a finite forbidden family
once the codegree density of the coloured class is known.
\begin{lemma}\label{lem:cyclic-local-rank}
Suppose a finitely described positive coloured model has pure facets completed in cyclic successor windows of size \(b\) on \(p\) colours, \(\gcd(b,p)=1\). Above \(b/p-\epsilon_0\), its exact codegree is realized by finitely many nonmembers, for some \(\epsilon_0(p,b)>0\). In rank three, if the inverse window entries have magnitude at most one, \(\epsilon_0=1/(100p^2)\) suffices.
\end{lemma}
\begin{proof}
Let \(h\) be the uniformity. In a host of codegree at least \((b/p-\epsilon)n\), hypergeometric concentration gives an induced subhypergraph on a sufficiently large constant number \(m\) of vertices with codegree ratio at least \(b/p-2\epsilon\). Fix these vertices as the anchor. Here a facet is an \((h-1)\)-set. For any tentative valid colouring, choose a facet in the nearest preceding class large enough to contain one. Its neighbours show that each window contains at least \((b/p-2\epsilon)m-p(h-2)\) anchor vertices. The window masses sum to \(b\). Their total absolute deviation from \(b/p\) is thus at most \(2p(2\epsilon+p(h-2)/m)\); invert the window matrix to see that every anchor class has a pure facet. In rank three use the entrywise inverse bound.

Choose one monochromatic \((h-1)\)-set of each colour in the anchor. The neighbourhoods containing a given vertex must form a subset of one of \(p\) distinct incidence patterns of size \(b\). A violation can be witnessed on boundedly many vertices. If a vertex lies in exactly \(b\) of these neighbourhoods, its colour is uniquely determined. Summing neighbourhood sizes leaves at most \(p\epsilon n\) unpinned. Provisionally assign compatible colours to them: every window has mass \(\ge b/p-\epsilon\), hence \(\le b/p+(p-1)\epsilon\). Inversion gives each colour mass \(1/p+O_{p,b}(\epsilon)\), with error at most \(2p\epsilon\) under the entrywise bound. Thus any two distinct windows' pinned intersection plus the unpinned set has mass \(<b/p-\epsilon\) for small \(\epsilon\) (under that bound it is at most \(b/p-1/p+(5p-1)\epsilon\)). So the pinned neighbours of an unpinned facet permit at most one pure colour. This is certified by one neighbour of each used colour and its pin. Two valid assignments on the same set differ at most \(p(p-1)(h-2)\) times, because \(h-1\) moves from a fixed colour to a different fixed colour would require two labels for their facet. Apply Lemma~\ref{lem:bounded-colour}.
\end{proof}

\begin{proof}[Proof of the finite-family assertion of Theorem~\ref{cg:main} for \(k=3\)]
The inverse of \(C\) is circulant with row and column sums \(1/b\).
Its entries take the two values \(\kappa\) and \(\kappa-1\), where
\(\kappa=(c+1/b)/p\). Indeed, the differences of the first row satisfy
\(d_j-d_{j-b}=\mathbf1_{j=p-1}-\mathbf1_{j=0}\). The smaller entries
occur at \(j=0,b,\ldots,(c-1)b\). On a bounded interval of class
indices only the diagonal entries have this value for large \(p\).

Choose distinguished classes \(T_i,A,B\) in a bounded interval. Give
each \(T_i\) its own block of auxiliary classes, and put \(A\) and then
\(B\) after all these blocks. Consider one pair with a vertex in \(A\)
and a vertex in \(T_i\). Begin with its window \(W_A\), and choose starts
\(P<T_i<D<A\) in block \(i\). Adding
\({\bf1}_{W_P}-{\bf1}_{W_D}\) inserts the interval
\(\{P+1,\ldots,D\}\) and removes its translate by \(b\). Further
differences \({\bf1}_{W_{P'}}-{\bf1}_{W_{D'}}\) with consecutive starts
may be inserted on disjoint intervals after \(D\) in the same block.
For large \(p\) the interval containing the starts is shorter than both
\(b\) and \(p-b\). The inserted intervals then lie outside \(W_A\) and
their translates lie inside \(W_A\). The resulting indicator \(w_i'\) thus
still has \(b\) entries equal to one. Its support contains \(T_i\) and
\(B\) and excludes \(A\) and all other \(T_j\). This agrees with the
pure-pair rule on repeated classes and keeps the rules for different
mixed pairs compatible. Since each window is a row of \(C\), the row
\(\eta_i=w_i'C^{-1}\) has coefficient \(+1\) at \(A\), and all its other
nonzero entries are signed coefficients in block \(i\).

Fix a large even integer \(J_0\) and distinct integers \(J_i\ge J_0\).
For each \(z\in B\), forbid \(K_{J_0+1}\) in \(L_H(z)[A]\) and
\(K_{J_i+1}\) in \(L_H(z)[A\cup T_i]\) for every \(i\). Write
\(s_i=1/J_0-1/J_i\). To prescribe a rational loss \(q_i>0\) in block
\(i\), write it as a finite sum of unit fractions, allowing repetitions.
For each summand \(1/M\) use a distinct auxiliary class with negative
coefficient in \(\eta_i\) and forbid \(K_{M+1}\) in the link of \(z\)
on that class. The graph Tur\'an bound gives missing-pair proportion at
least \(1/M-o(1)\), attained asymptotically by a balanced complete
\(M\)-partite link. The additional pairs of starts above provide one
negative coefficient per summand.

Choose rational weights \(\theta,w_i<1\) close to one. For a designated
class \(H\) write \(\beta_A=\theta\) and \(\beta_{T_i}=w_i\). In links
rooted at vertices of \(H\), impose clique restrictions of total loss
\(1-\beta_H\) on preceding auxiliary classes outside all the supports
of the \(\eta_i\). Set \(\rho=\sum_H(1-\beta_H)\). Let \(R\) have one
column \(r_H\) per designated class, giving its losses on these
auxiliary rows, and let the corresponding column of \(E\) be the
coordinate vector \(e_H\). The matrix \(D=C-RE^{\mathsf T}\) incorporates
these losses. Using \(E^{\mathsf T}C^{-1}r_H=\kappa(\sum r_H){\bf1}\)
in the rank update gives, for \(\omega={\bf1}^{\mathsf T}D^{-1}\) and
\(v=e_B^{\mathsf T}D^{-1}\),
\begin{equation}
\begin{gathered}
\omega_j=\alpha_j/[b(1-\kappa\rho)],\quad v_j-b\kappa\omega_j=(C^{-1})_{B,j}-\kappa,\\
{\bf1}^{\mathsf T}D^{-1}{\bf1}=p/b+\rho/[b^2(1-\kappa\rho)],\quad
(D^{-1}{\bf1})_B=1/[b(1-\kappa\rho)],\quad w_i'D^{-1}=\eta_i .
\end{gathered}\label{eq:cgproof-2}
\end{equation}
Here \(\alpha_H=\beta_H\), all other local \(\alpha_j=1\), and globally \(\alpha_j\ge1-\rho\).

For the upper bound it suffices to consider \(n\)-vertex hosts of normalized codegree at least \(b/p-1/(100p^2)\). Each class then has \(\Omega(n)\) vertices. Indeed, for each window \(W_j\), use a pair from the nearest preceding class with at least two vertices. Passing to \(W_j\) loses at most \(p\) intervening vertices. The window proportions sum to \(b\), and their total absolute deviation from \(b/p\) is at most \(1/(50p)+o(1)\). Applying \(C^{-1}\) gives the claim about the individual classes.

In the graph links of vertices of \(B\), write the average proportions of
missing pairs within \(A\) and within \(T_i\) as \(1/J_0+a_0\) and
\(1/J_i+a_i\), respectively, with \(a_j\ge0\). Write \(S_i\) for the
average missing proportion between \(A\) and \(T_i\). Give the two classes total weights \(u,v\ge0\), respectively.
The weighted graph Tur\'an bound on \(A\cup T_i\), averaged over the
roots in \(B\), gives
\[
  (a_0+s_i)u^2+2(S_i-1/J_i)uv+a_iv^2\ge0.
\]
Varying \(u/v\) gives \(S_i\ge1/J_i-\sqrt{(a_0+s_i)a_i}\),
with the same bound by continuity if a diagonal coefficient is zero. Counting diagonals as missing
costs \(O(1/n)\). Denote the normalized codegree after this error by \(d\).

The pure-pair inequalities now read \(Dx=d{\bf1}+x_Bz\). Here \(z\)
records the within-class losses on \(A\) and the \(T_i\), the prescribed
auxiliary losses totalling \(q_i\) in block \(i\), and a nonnegative
slack vector \(f\). For the mixed pair from \(A\) and \(T_i\) the
corresponding inequality is \(w_i'x\ge d+x_B S_i\). Substituting for
\(x\) and using \eqref{eq:cgproof-2} gives
\begin{equation}
q_i\le a_0+s_i+\sqrt{(a_0+s_i)a_i}+\eta_i f. \label{eq:cgproof-3}
\end{equation}
To see which losses attain equality, first set \(f=0\) and \(a_0=t>0\).
In the range \(0<t+s_i<q_i\), the smallest admissible \(a_i\) in
\eqref{eq:cgproof-3} is \((q_i-t-s_i)^2/(t+s_i)\). Thus the weighted
excess loss \(\theta a_0+\sum_i w_i a_i\) to be minimized is
\[
 \begin{aligned}
 V(t)&=\theta t+\sum_iw_i(q_i-t-s_i)^2/(t+s_i)\\
     &=c_0t+\sum_iw_iq_i^2/(t+s_i)+\sum_iw_i(s_i-2q_i),
 \qquad c_0=\theta+\sum_iw_i .
 \end{aligned}
\]
Differentiation gives \(V'(t)=c_0-F(t)\), where
\(F(t)=\sum_iw_iq_i^2/(t+s_i)^2\). At the admissible positive root
of \(F(t)=c_0\), put \(\lambda_i=q_i/(t+s_i)-1\). The minimizing
choice above is \(a_i=\lambda_i^2(t+s_i)\). The positive multipliers
\(\mu_i=2w_i\lambda_i\) match the derivatives since
\(\theta=\sum_i w_i(2\lambda_i+\lambda_i^2)\), and
\(m_0=\sum_i\mu_i<\theta\). Concavity of the geometric means in
\eqref{eq:cgproof-3} now bounds the weighted excess loss with slack
below by \(V(t)-m_0\sum_j f_j\). Put
\(L_*'=\theta/J_0+\sum_iw_i/J_i+\sum_iq_i+V(t)\), \(L'=\sum_j\alpha_j z_j\), \(s=\sum f_j\). Then \(L'\ge L_*'+(1-\rho-m_0)s\) and \(vz\ge\kappa L'/(1-\kappa\rho)-s\).

Solving in coordinate \(B\) gives
\(1/d={\bf1}^{\mathsf T}D^{-1}{\bf1}+(D^{-1}{\bf1})_B\,\omega z/(1-vz)\).
The denominator is positive, and if \(L_*'<1-\rho-m_0\) this is at least the value for \(L'=L_*',s=0\) (cross multiplication gives \(L'-L_*'\ge L_*'s\)).

To build links attaining the upper bound, give the vertices of \(A\) one of \(J_0\) auxiliary colours (shades) with probabilities \(P_j=(1\pm\sqrt{J_0t})/J_0\) using each sign \(J_0/2\) times. Extend by zeros to \(J_i\) and give \(T_i\) probabilities \(Q_{ij}=1/J_i-\lambda_i(P_j-1/J_i)\). For \(j\le J_0\), \(J_iP_j=(J_i/J_0)(1\pm\sqrt{J_0t})\) is close to one in the region chosen below, so \(J_iQ_{ij}=1-\lambda_i(J_iP_j-1)>0\). For \(J_0<j\le J_i\), \(P_j=0\) and \(Q_{ij}=(1+\lambda_i)/J_i>0\). The probabilities of equal shades within a class and across the two classes attain equality in \eqref{eq:cgproof-3}. For each root form the multipartite link by joining vertices of different shades. For this root use the same shades on \(A\) in every \(A\cup T_i\), and independent uniform shades for the other clique restrictions. Make the choices independently for different roots.

It remains to check the class proportions. For the solution with no slack, put \(L=\rho+L_*'\). Since \(|(C^{-1})_{ij}|\le1\) and the tilt columns have total mass \(\rho\), the rank update used in \eqref{eq:cgproof-2} gives
\[
 \|bD^{-1}{\bf1}-{\bf1}\|_\infty\le\frac{\rho}{1-\kappa\rho},
 \qquad \max_{i,j}|(D^{-1})_{ij}|\le1+\frac{\rho}{1-\kappa\rho}<2.
\]
At equality \(x_B=d/[b(1-\kappa L)]\) and \(\|z\|_1\le L_*'/(1-\rho)\). Thus \(x=dD^{-1}{\bf1}+x_BD^{-1}z\) satisfies
\[
 \left\|(b/d)x-{\bf1}\right\|_\infty
 \le \frac{\rho}{1-\kappa\rho}
      +\frac{2L_*'}{(1-\rho)(1-\kappa L)}<0.02
\]
in the region \(L<0.004\) used below. The normalized masses are therefore positive and of order \(1/p\). Concentration gives the required number of neighbours simultaneously for pairs within a class and for the specified mixed pairs. Every other mixed pair can be completed in all but boundedly many classes. Together with the upper bound this gives the exact normalized codegree
\begin{equation}
d=\left(\frac pb+\frac{L}{b^2(1-\kappa L)}\right)^{-1}
=\frac bp-\frac{L}{p^2(1-cL/p)},\qquad L=\rho+L_*'. \label{eq:cgproof-4}
\end{equation}
For \(L<.004\) and sufficiently large \(p\), this value lies within \(1/(100p^2)\) of \(b/p\). Since \(\gcd(p,b)=1\), Lemma~\ref{lem:cyclic-local-rank} now gives a finite forbidden family at this density.

To place the conjugates of the density, we first arrange the number of real roots of \(F(t)=c_0\). For each \(1\le h\le r\) choose unit prototype weights and \(0<\lambda_i<1/2\) with \(\sum_i(2\lambda_i+\lambda_i^2)=1\), taking \(\lambda_1\in(.35,.42)\) and the rest small so the gap \(1-m_0=\sum\lambda_i^2>.12\). Put \(q_i=(1+\lambda_i)(t_0+s_i)\). Arrange the \(s_i\)'s in \(h\) tight consecutive clusters at widely separated positive scales below \(.01/J_0\), with \(J_0t_0\ll .01\). These scales are approximated by the mesh \(J_0s_i=1-J_0/J_i\) for large even \(J_0\). Then \(F(t_0)=r+1\). Between consecutive poles \(F\) is strictly convex with infinite endpoints. Inside a cluster a pole term exceeds \(r+1\) across the gap. Between clusters, at \(-\sqrt{(t_0+s_i)s_{i+1}}\) the terms \(j\le i\) tend to zero as the separating ratio grows, while the others approach their values at \(t_0\), so the value is below \(r+1\). Each outer ray has one root by monotonicity. Hence there are \(2h\) simple real roots. The positive-root conditions, positivity of the \(Q_{ij}\), and \(L<.004,\ L_*'<1-\rho-m_0\) persist: the untilted losses are \(O_r(1/J_0)\) and the multiplier gap is bounded below. Take rational \(q_i\) close to their prototypes and move the weights slightly below one, obtaining an open parameter region.

Write \(d_K=[K : \Q]\). With the shifts \(s_i\) and losses \(q_i\)
fixed, we vary the residues \(r_i=w_iq_i^2\). Our first requirement
is that the chosen conjugates of \(t\) all satisfy \(F(t)=c_0\).
At these stationary points, the values of \(V\) determine the losses
in \eqref{eq:cgproof-4}. For \(r\ge2d_K-1\), distinct points \(\tau_j\) off the poles give rank \(d_K\) for \(((\tau_j+s_i)^{-2})_{j,i}\): a dependence produces \(\sum_j a_j/(\tau_j+z)^2\) with numerator of degree \(\le2d_K-2\) vanishing at the \(r\) distinct \(s_i\). For distinct real roots \(t,u\), use the stationary equation in the difference \(V(t)-V(u)\) to obtain
\begin{equation}
V(t)-V(u)=\frac{t-u}{2}\sum_i r_i\big((t+s_i)^{-1}-(u+s_i)^{-1}\big)^2.
\label{eq:cgproof-6}
\end{equation}
For a nonreal root \(t=x+iy\), the real part of \(F(t)=c_0\) gives
\begin{equation}
\Im V(t)=-2y^3\sum_i r_i/((x+s_i)^2+y^2)^2\ne0. \label{eq:cgproof-7}
\end{equation}

For \(d_K\ge2\) choose \(r\ge2d_K-1\) and a root count allowing the required real and nonreal embeddings (use the least even count at least the real signature). Choose a rational prototype with all roots simple, avoiding the finitely many critical values of \(F\). Select a conjugation-stable set of \(d_K\) roots with the admissible positive root distinguished, then approximate them by the conjugates of \(t\in K\) using Lemma~\ref{lem:choose-generator}. The interpolation rank remains \(d_K\), so \(d_K\) of the \((t+s_i)^{-2}\) are a rational basis. Fix \(c_0\) and the other residues, and solve \(\sum_i r_i(t+s_i)^{-2}=c_0\) in that basis. The residues are rational and close to their prototypes, hence feasible. The rightmost root \(t\), used in the construction, gives a value \(V(t)\) distinct from every other conjugate by \eqref{eq:cgproof-6}--\eqref{eq:cgproof-7}; thus \(\Q(V(t))=K\). Equation \eqref{eq:cgproof-4} preserves this field and reverses the real ordering of the losses. Thus the density is the least of its real conjugates, and all its conjugates tend to \(A_0/c\) as \(p\to\infty\). For \(K=\Q\) use the uncapped cyclic model of density \(b/p\), by summing the pure inequalities and balancing.
Choose \(A_0/c\) within \(\varepsilon/2\) of the prescribed \(u\), and
then take \(p\) sufficiently large for all embeddings of \(K\).
\end{proof}

For the minimum pair-degree construction in Section~\ref{ch:quartic} we shall need finer control of the second-order term in this codegree.
\begin{theorem}\label{ch:region}
For each \(K\) there is a bounded nonempty open \(U\subset K_\infty\) whose elements in \(K\) are generators of positive distinguished value. For each fixed reduced rational \(A_0/c\in(0,1)\) there is a constant \(C=C(K,A_0/c)\) with the following property. If \(g\in K\cap U\) and a positive integer \(p\) satisfies \(A_0p\equiv1\pmod c\) and \(p\ge C(1+\log H(g))\), then \(b/p-g/p^2\), with \(b=(A_0p-1)/c\), is an exact finite-family codegree density. Here \(H\) is coordinate height in a fixed rational basis.
\end{theorem}

\begin{proof}
Write \(d_K=[K : \Q]\) and take \(r\ge2d_K\). Fix the admissible root
\(t\in K\) and \(c_0\) from the preceding construction (for \(\Q\)
there is just one stationarity condition). We now choose the residues
to prescribe the loss at this root as well. A small rational loss \(z>0\) at extra pure rows rooted at \(B\) outside all mixed/tilt supports adds \(z\) to \(L\). Holding stationarity, the targets are
\[
\sum_i r_i(t+s_i)^{-2}=c_0,\qquad
L=C+c_0t+\sum_i r_i((t+s_i)^{-1}-2/q_i)+z,\quad
C=r+1-c_0+c_0/J_0+\sum_iq_i.
\]
These two maps from the rational parameters onto \(K\oplus K\) are jointly surjective: use the \(z\) column to remove constants \(-2/q_i\); a complex row dependence gives \(\sum_j(a_j/(X+t_j)^2+b_j/(X+t_j))\) of numerator degree \(\le2d_K-1\), vanishing at the \(r\) distinct \(s_i\). Here \(t_j\) are the conjugates of \(t\). A fixed rational right inverse therefore realizes all nearby losses in \(K\) with the stationary equation fixed. Start at a small positive rational \(z\); all inequalities are strict and the distinguished loss is distinct from its other conjugates, an open condition.

Choose a nonempty open set \(U\) with closure compactly contained in this region. For the prescribed \(g\) and \(p\), take \(L_p=g/(1+cg/p)\); \eqref{eq:cgproof-4} becomes \(b/p-g/p^2\). Uniform convergence at all embeddings preserves admissibility for large \(p\). Inverting the multiplication matrix for \(1+cg/p\) by adjugates gives heights \((pH(g))^{O(1)}\), as does the fixed right inverse. A positive rational \(a/e<1\) uses \(O(1+\log e)\) unit fractions: take \(j=\lceil\log_2e\rceil\), \(2^ja=ke+r\), and expand \(k/2^j+r/(2^je)\) in binary numerators. The auxiliary interval therefore has length \(O(1+\log p+\log H(g))\). For fixed \(A_0/c\), both \(b\) and \(p-b\) grow linearly with \(p\). Taking \(C=C(K,A_0/c)\) large ensures that the interval is shorter than both and that \(c\) times its length is below \(p-c\). Thus the required off-diagonal high inverse entries and localization margin persist.
\end{proof}

We now pass from \(r\)-graphs to \((r+1)\)-graphs by placing an \(r\)-graph from the preceding construction in a vertex link.
\begin{theorem}\label{cg:thm:next_rank}
For a finite-family \(r\)-graph codegree density \(\beta<1\), \(r\ge3\), fix coprime \(1\le b<p\), \(p\ge4\). For every sufficiently large integer \(t\), put \(c=(1-\beta)^{-1}\), \(x_0=(t+c)^{-1}\), \(x_i=(t+1)^{-1}\) otherwise. If \(C\) is the cyclic matrix of windows \(\{i,\ldots,i+b-1\}\) and \(R=C-\mathrm{diag}(x_i)\), then \(\alpha_t=({\bf1}^{\mathsf T}R^{-1}{\bf1})^{-1}\) is an exact rank-\((r+1)\) finite-family codegree density, tends to \(b/p\), and generates \(\Q(\beta)\).
\end{theorem}

\begin{proof}
Index the vertex classes \(V_i\) cyclically. An \(r\)-set contained in \(V_i\) has neighbours only in \(V_{i+1},\ldots,\allowbreak V_{i+b}\). For \(z\in V_{i+1}\), require \(L_H(z)[V_i]\) to avoid the \(t\)-fold codegree suspension of the source family if \(i=0\), and the corresponding suspension of a single edge otherwise. By Lemma~\ref{lem:family-suspension}, the suspended family has codegree
density \(1-x_i\). To attain this density in the link, put \(c_i=c\) for
\(i=0\) and \(c_i=1\) otherwise. Give each of the \(t\) independent
classes proportion \(1/(t+c_i)\) and the core proportion \(c_i/(t+c_i)\).
Choose a near-extremal core avoiding the input family and include all
crossing edges. This link avoids the suspended family, and its missing
edges lie within individual parts. Its ordinary edge-density deficit
is therefore at most \((t+c_i^r)/(t+c_i)^r+o(1)=O_{r,c}(t^{1-r})+o(1)\).

A critical facet of type \((r-1)V_i+V_{i+1}\) completes only in \(\{i,\ldots,i+b-1\}\); otherwise include all allowed edges. Critical prescriptions are compatible: the only two-part ambiguity \(r=3,\,2+2\) would require opposite cyclic arcs.

When the normalized codegree is close to \(b/p\), the window argument for \(r\)-sets in one class gives \(|V_i|=\Omega(n)\) for every \(i\). A minimum facet in each capped link gives \(Ra\ge d{\bf1}-o(1)\) for proportions \(a\) and normalized degree \(d\). The window eigenvalues are \(b\) at \(1\) and \((1-\xi^b)/(1-\xi)\ne0\) otherwise. Both inverse-sum vectors of \(R\) are positive for large \(t\), giving the upper bound \(\alpha_t\).

For a matching construction take proportions \(R^{-1}{\bf1}/({\bf1}^{\mathsf T}R^{-1}{\bf1})=p^{-1}{\bf1}+O(t^{-2})\); \(\alpha_t=b/p-1/(pt)+O(t^{-2})\). Independent relabellings of these extremal link constructions give the critical degrees and, by concentration, pure degrees at least \(b/p-O(t^{1-r})-O(t^{-2})-o(1)\).

Other facets have at least \(b\) unrestricted completion parts. For two-part facets only an own-class completion creating a critical profile can be lost, except for the profile \((r-1)V_i+V_j\) with \(j\ne i+1\), which can also lose \(V_i\) by the pure rule and \(V_{i+1}\) by the critical rule. A loss of \(V_j\) in this last case requires \(r=3,b=1\). For \(b=p-1\) the pure rule excludes no distinct colour, leaving at most one loss. For three-part facets only existing parts may be lost, at most one if \(r\ge4\), and at most two if \(r=3\) (arcs of the directed cycle on three colours). When \(b=p-1\), a three-part loss additionally requires the predecessor outside the window, so only the middle vertex of a directed length-two path can be lost. Four or more parts cause no exclusion. This proves the matching lower bound. Lemma~\ref{lem:cyclic-local-rank} gives finite witnesses in the new uniformity, since the window length \(b\) is coprime to the number \(p\) of classes.

It remains to trace the generated field through this construction. For \(x_*=(t+1)^{-1}\), \(z=((C-x_*I)^{-1})_{00}\), \(u=x_0-x_*\), a rank-one update gives
\[
 1/\alpha_t=p/(b-x_*)+u/((b-x_*)^2(1-uz)).
\]
Put \(\Delta=b-x_*\) and \(z'=z-1/(p\Delta)\). Substituting \(u=-\beta/((t+1)(t+1-t\beta))\) and clearing denominators gives
\[
 \alpha_t=\frac{\Delta}{p}\,
 \frac{(t+1)^2-[t(t+1)-z]\beta}
      {(t+1)^2-[t(t+1)-z']\beta}.
\]
Since \(z-z'=1/(p\Delta)>0\), this rational fractional-linear map is invertible and has positive derivative on each real interval avoiding its pole. In particular \(\Q(\alpha_t)=\Q(\beta)\). As \(C\) is invertible, \(z'=O_{p,b}(1)\), and the pole is \(1+1/t+O_{p,b}(t^{-2})\). For Theorem~\ref{cg:main} we place all conjugates of the input sufficiently close to the chosen target in \((0,1)\) that its real conjugates lie in \((0,1)\). For large \(t\) this interval avoids the pole, so their order is preserved. At each embedding \(u\to0\), hence all conjugates of \(\alpha_t\) tend to \(b/p\). This gives both simultaneous placement and the least-real-conjugate property in the next uniformity.
\end{proof}

In the passage to a single forbidden hypergraph we also need a finite-family criterion for classes defined by weak colourings.
\begin{theorem}\label{weakloc:thm:loc}
A class \(\D\) represented by \(p\ge2\) independent parts and finitely many positive coloured prohibitions has finite exact codegree witnesses if its codegree density \(\beta>1-1/(p-1/2)\).
\end{theorem}
\begin{proof}
Fix a codegree ratio strictly between the stated threshold and \(\beta\), and write \(h\) for the uniformity. Suspending one edge \(p-2\) times~\cite{ADLY} gives an \(h\)-graph of weak chromatic number \(p\) and codegree density at most \(1-1/(p-1)\), since deleting any independent set leaves the input at each step. The host thus contains \(h-1\) disjoint copies. In any tentative colouring each copy uses all \(p\) colours, so we can fix a monochromatic \((h-1)\)-set \(P_c\) of each colour \(c\) in this anchor.

Put \(\mu=1-\delta_{\co}(H)/n\) and \(T_c=V(H)\setminus N_H(P_c)\). Every vertex of colour \(c\) belongs to \(T_c\). Unless a bounded incidence obstruction already rules out the colouring, these sets cover \(V(H)\). The part of \(T_c\) lying in no other \(T_j\) has at least \((\mu-(p\mu-1))n>\mu n/2\) vertices. Indeed, the unused allowance \(p\mu n-\sum|T_c|\) and the total overlap \(\sum|T_c|-n\) sum to \((p\mu-1)n\). A vertex in this private region is forced to have colour \(c\), as one incident edge rules out each other colour. The neighbourhood of any \((h-1)\)-set meets all but at most one private region. Hence an \((h-1)\)-set of unpinned vertices has at most one possible colour if it is monochromatic. One pinned neighbour of each colour met certifies this restriction. The bound on the distance between valid colourings now gives the claim by Lemma~\ref{lem:bounded-colour}.
\end{proof}

\subsection{A single forbidden hypergraph}\label{sec:codegree-single}\label{cg:part:singleton}
Throughout this subsection \(h\) denotes the uniformity. Let \(\D\) be a hereditary class represented by \(p\) ordered independent parts and finitely many forbidden coloured configurations. To pass to a single forbidden \(h\)-graph, we take a product whose factors each consist of a base hypergraph of high minimum codegree and a collection of attached hypergraphs. The upper bound embeds a factor. For the lower bound we show that a colouring of the product would induce a forbidden colouring on an attachment.
\begin{theorem}[Singleton conversion]\label{cg:thm:conversion}
If the codegree of \(\D\) is \(\beta\) and a finite family \(\A_0\) outside it has exact density \(\beta\), then for all large \(q\) one graph has codegree density \(\alpha_q=1-(q+(1-\beta)^{-1})^{-1}\), generating \(\Q(\beta)\).
\end{theorem}

We begin with the independent sets in the product of the \(m\)-vertex base hypergraphs. Write \(d=(1-\mu)m\) for the minimum codegree bound.

\begin{lemma}\label{cg:lem:ind}
For small \(\mu\), a product independent set has mass \(\alpha\le\mu+10\mu^2\).
\end{lemma}
\begin{proof}
First consider a graph factor. We choose weights on its edges so that all vertices have the same weighted degree. Put \(\nu=\mu+10\mu^2,\ R=(1-\nu)m\). Remove a symmetric fractional subgraph of degrees \(s_i=\deg(i)-R\in[10\mu^2m,\nu m]\). The bipartite double-cover cut condition is \(s(X)+s(Y)-\sum_i s_i\le e(X,Y)\): if \(|Y|\ge(\mu+\nu)m\) then \(e(X,Y)\ge\nu m|X|\), and symmetrically for \(X\); otherwise \(s(X)+s(Y)\le2(\mu+\nu)\nu m^2<\sum_i s_i\). Symmetrize the flow. Put \(\tau=\nu/(1-\nu)\) and let \(B\) be the remaining normalized kernel. The symmetric nonnegative matrix \(J-RB\) has row sum \(\nu m\), so its operator norm is at most \(\nu m\). On the orthogonal complement of the constants it equals \(-RB\). Thus the centred norm of \(B\) is at most \(\tau\).

For graphs the product kernel vanishes on the set, so \(0\ge\alpha^2-\tau\alpha(1-\alpha)\). For the induction step, fix a vertex \(z\) of the independent set and take the product of the links at its coordinates, on \(Q_z=\prod_i(V_i\setminus\{z_i\})\), with defect \(\mu'=1-d/(m-1)\). The restriction of the independent set to this product has mass \(\le\mu'+10(\mu')^2\). The tensor \(P\) of normalized complete-graph operators has least eigenvalue \(-1/(m-1)\); combining \(\langle f,Pf\rangle\le\alpha(\mu'+10(\mu')^2)\) with its spectral lower bound yields the claim.
\end{proof}
\begin{lemma}\label{cg:lem:tail}
For fixed uniformity \(h\), \(\tail(f)=O_h((\mu-\alpha)_++\mu^2+1/(m\sqrt\alpha))\) for a nonempty independent set, uniformly in the number of factors.
\end{lemma}
\begin{proof}
The case \(\alpha\le\mu/2\) is trivial. In the graph case, let \(\nu,\tau\) be as in the proof of Lemma~\ref{cg:lem:ind}. The spectral bound there improves by \((\tau-\tau^2)\tail(f)\), giving \(\tail(f)\le\alpha(\nu-\alpha)/[\nu(1-\tau)]\). Induct as above: \(\alpha_z=Pf(z)\le\mu'+10(\mu')^2\), \(\E(\alpha_z\mid f(z)=1)\ge\alpha-1/(m-1)\). On slices of mass below \(\mu'/2\) use the trivial bound, on the rest induction. The resulting mean slice tail is \(O_h((\mu-\alpha)_++\mu^2+1/(m\sqrt\mu))\).
Use the same four-block statistic as in the ordinary conversion, now drawing \(X,Y\) independently uniform on \(Q_z\). Its conditional expectation is at most the slice tail. The centred pair-to-root channel in one coordinate has norm \(\le2/(m-2)\): the conditional mean of a centred vertex function is the negative sum at the one or two excluded symbols, divided by the number available. Tensorization gives expectation change \(O(1/(m\sqrt\alpha))\) on unconditioning the root. The unconditional multiplier is \(1+3\rho^j-4(\rho+(1-\rho)/4)^j\), \(\rho=(m-1)^{-2}\), zero for \(j=0,1\) and at least \(1/2\) otherwise. This transfers the estimate.
\end{proof}

We apply the tail estimate to the colour classes of a weak colouring to find a coordinate that predicts the colours of most vertices.
\begin{lemma}\label{cg:lem:colors}
Fix \(h,p\), let \(\mu>0\) be sufficiently small, and take \(m\) sufficiently large in terms of \(h,p,\mu\), in particular \(m\mu^{5/2}\ge1\). If a weak proper \(t\)-colouring has \(t\mu-1\le p\mu\), a modal classifier on one coordinate has error \(O_{h,p}(\mu)\). Except for \(O_{h,p}(1)\) colours \(c\), it has disjoint slices \(S_c\) of size at least \(\mu m/3\), each row seeing \(c\) with probability at least \(1-1/(2(h-1))\).
\end{lemma}
\begin{proof}
Since \(\alpha_c\le\mu+10\mu^2\), the sum of mass deficits is \(O_p(\mu)\). For \(\alpha_c<\mu/2\), the tail of the colour indicator \(f_c\) is at most \(\alpha_c\le\mu-\alpha_c\). For the remaining colours sum Lemma~\ref{cg:lem:tail}. Since \(t\le1/\mu+p\), their \(\mu^2\) terms sum to \(O_p(\mu)\) and their finite-order errors sum to \(O_p(1/(m\mu^{3/2}))\). Hence
\[
 \sum_c\tail(f_c)=O_{h,p}\!\left(\mu+\frac{1}{m\mu^{3/2}}\right)=O_{h,p}(\mu).
\]
By Lemma~\ref{lem:one-coordinate}, a colouring determined by one coordinate has the asserted error. Summed over the colours, the measures of the symmetric differences between the true colour classes and their predicted fibres equal twice this error. Discard colours of mass \(<\mu/2\) or symmetric difference \(>\mu/(100h)\). Only boundedly many are discarded. In each remaining colour remove the rows on which its conditional probability is too small. At least \(\mu m/3\) rows remain.
\end{proof}

It remains to describe the witnesses that will be attached to the factors.
Write \(\D^{+r}\) for the class obtained by joining a core in \(\D\) to
\(r\) independent classes, and put \(c_*=(1-\beta)^{-1}\).
\begin{lemma}[Adding independent classes]\label{lem:joined-codegree-witnesses}
Under the hypotheses of Theorem~\ref{cg:thm:conversion}, the extremal
codegree density of \(\D^{+r}\) is \(\alpha_r=1-1/(r+c_*)\) for each
\(r\ge0\). There is a finite family \(\A_r\) outside \(\D^{+r}\) with
\(\pi_{\co}(\A_r)=\alpha_r\).
\end{lemma}

\begin{proof}
Giving the core proportion \(c_* /(r+c_*)\) and each independent class proportion \(1/(r+c_*)\) gives codegree density \(\alpha_r=1-1/(r+c_*)\). Conversely, at minimum codegree \(d\) an independent class has at most \(n-d+O(h)\) vertices, and the codegree induced on the core is at least its order minus \(n-d\).

At the next step apply Lemma~\ref{lem:family-suspension} with
\(\ell=h-1\) to \(\A_r\). Take \(h-1\) root members and attach a member at every \((h-1)\)-set consisting of one vertex from each root. Deleting any independent set leaves an input member. For the upper bound at the next step, the neighbourhood of each such \((h-1)\)-set in a host above \(\alpha_{r+1}\) has induced codegree ratio above \(2-1/\alpha_{r+1}=\alpha_r\), even after boundedly many vertices are deleted.
Starting with the given family \(\A_0\), this gives the upper bound
at each step. A copy of a witness in a member of \(\D^{+(r+1)}\) would
leave a witness from \(\A_r\) inside a member of \(\D^{+r}\) on deleting
the last independent class, contrary to the induction hypothesis.
The balanced construction therefore gives the matching lower bound.
\end{proof}

\begin{proof}[Proof of Theorem~\ref{cg:thm:conversion}]
Use the families supplied by Lemma~\ref{lem:joined-codegree-witnesses}.
Let \(L_0\ge h\) bound the coloured witness orders and \(K_0\) bound bad colours in Lemma~\ref{cg:lem:colors}. Fix \(L\ge(p+K_0)L_0\), large \(q>K_0\), and \(1/(q+c_*)<\mu_0<1/(q+1/2)\) (here \(1\le c_*\le p\)). Take \(m\) a sufficiently large multiple of \(L\) and set \(d=\lceil(1-\mu_0)m\rceil\). Then \(\mu=1-d/m\to\mu_0>0\) as \(m\to\infty\), so in particular \(m\mu^{5/2}\ge1\) as required in Lemma~\ref{cg:lem:colors}. Take all \(m\)-vertex base hypergraphs of minimum codegree at least \(d\). For every sequence \(J\) of \(i\le q\) sets of size \(h-1\) in a base hypergraph, choose a member \(A_J\in\A_{q-i}\) and join each of its vertices to every facet of \(J\). Use disjoint copies for different sequences. Include every choice of attachments and take the product \(P\) of the resulting hypergraphs.

For the upper bound, use concentration to select an \(m\)-vertex base hypergraph of minimum codegree at least \(d\) in a host \(H\) whose codegree ratio exceeds \(\alpha_q\) by a fixed positive margin. The common neighbourhood of \(i\) such \((h-1)\)-sets, apart from used vertices, has order at least \(n-i(n-\delta_{\co}(H))-O(1)\), and minimum codegree at least its order minus \(n-\delta_{\co}(H)\), giving ratio above \(\alpha_{q-i}\). Embed all attachments greedily. The fixed finite family of factors thus has density \(\le\alpha_q\); Lemma~\ref{lem:product-upper} gives the product upper bound.

For the lower bound, suppose that \(P\) admits a representation in \(\D^{+q}\). The product of the base hypergraphs then has \(t=q+p\) weak colours and \(t\mu-1\le p\mu\) at \(\mu=1-d/m\).
Fix the coordinate given by Lemma~\ref{cg:lem:colors}. A row now corresponds to a vertex of this factor, and an outer tuple specifies the vertices in all the other factors. We use the sets \(S_c\) of rows to restrict the colours of neighbours of an \((h-1)\)-set.
For fixed \(L\mid m\), Lemma~\ref{lem:uniform-cliques} applies for small \(\mu\) (averaging missing codegrees bounds missing vertex degrees). On the outer product, the complete-product operator has centred norm \(1/(m-1)\). Chebyshev and a union bound over the at most \(m\) slice rows give a set \(U\) of measure \(1-O_h(1/m)\) where all smoothed slice-colour probabilities exceed \(1-1/(h-1)\).
Let \(T\) be an \((h-1)\)-set in \(S_c\). If \(T\cup\{x\}\) is an edge of the chosen factor, then the product vertex \((x,w)\) cannot have colour \(c\) for any \(w\in U\). This holds also when \(x\) is an attached vertex. Indeed each outer link rooted at \(w_i\) has an ordered edge law with uniform marginals on the remaining vertices: induct by choosing a uniform first vertex, down to a perfect matching in the double cover of a graph of minimum degree at least half its order. The product link laws and a union bound over the rows of \(T\) would otherwise give a monochromatic edge.

In this factor choose an \((h-1)\)-set of base vertices from the slice of each good pad colour, and let \(J\) be the resulting sequence. If \(k\le K_0\) pads are bad, the neighbourhood restriction just proved leaves at most \(k+p\) colours on the rows of \(A_J\in\A_k\) over \(U\). Align \(L\) complete positions in the outer factors using the uniform clique laws; all can lie in \(U\) since \(L\) is fixed and \(m\) large. Each attachment row's modal label occurs at least \(L_0\) times. Any monochromatic edge or forbidden coloured witness in these labels lifts using distinct outer positions. The resulting row colours would thus put \(A_J\) in \(\D^{+k}\), contrary to the choice of \(\A_k\) (with a sole factor use the empty-sequence attachment). Thus every member of \(\D^{+q}\) is \(P\)-free. The balanced examples in
this class give \(\pi_{h-1}(P)\ge\alpha_q\), as required. The fractional-linear transfer retains the field and preserves the real ordering when the input conjugates lie in \((0,1)\), which the preceding construction ensures.
\end{proof}

\begin{proof}[Proof of Theorem~\ref{cg:main}]
The three-uniform construction in Section~\ref{sec:codegree-classes}
gives the finite-family assertion for \(k=3\). To pass up one
uniformity, choose a reduced fraction \(b/p\) within \(\varepsilon/2\)
of \(u\), with \(p\ge4\). By Theorem~\ref{cg:thm:next_rank} and its
proof, the new density generates the same field, all its conjugates
approach \(b/p\) as \(t\) grows, and the order of the real conjugates
is preserved when they lie in \((0,1)\). Starting sufficiently close
to \(u\), we may thus iterate the construction up to uniformity \(k\).

Each resulting class has independent colour classes, finitely many
coloured prohibitions and a finite exact witness family. Apply
Theorem~\ref{cg:thm:conversion} to its density \(\beta\). The singleton
densities \(\alpha_q\) tend to one and generate \(K\). The transfer is
increasing on \((0,1)\) and is an invertible rational fractional-linear
map, so the least real conjugate remains the least real conjugate.
\end{proof}

The codegree theorem also allows us to specify the degree of \(K\), its number of real embeddings and the Galois group of its normal closure.

\begin{corollary}\label{cor:codegree-signature}
Let \(d\ge2\) and \(1\le s\le d\) satisfy \(s\equiv d\pmod2\), and let
\(L\) be a number field. The field \(K\) in Theorem~\ref{cg:main} may
be chosen to have degree \(d\) and exactly \(s\) real embeddings. Its
normal closure can at the same time be chosen to have Galois group
\(S_d\) over \(\Q\) and to intersect the normal closure of \(L\) only
in \(\Q\).
\end{corollary}

\begin{proof}
Prescribe squarefree reduction types \((d),(d-1,1),(2,1^{d-2})\) at distinct primes \(>d\) (only irreducibility is needed if \(d=2\)). At an auxiliary odd prime \(p_0>d\) unramified in the normal closure of \(L\), prescribe \((x^2-p_0)\prod_{j=1}^{d-2}(x-j)\) modulo \(p_0^2\); its discriminant has valuation one at \(p_0\). By CRT approximate the coefficients of \(M^dg(x/M)\) within bounded errors, where \(g\) is monic real with \(s\) simple real roots and the rest simple nonreal pairs. Rescaling back preserves signature for large \(M\). The modular types give transitivity, a \((d-1)\)-cycle and a transposition, hence \(S_d\). Its sign quadratic field ramifies at \(p_0\); every nontrivial Galois subextension contains it because every proper normal subgroup is contained in \(A_d\). The normal closures therefore have trivial intersection. Take any real root as the embedding and apply Theorem~\ref{cg:main}.
\end{proof}

\section{Intermediate minimum degrees}\label{sec:intermediate-cases}
Write \(s=k-\ell\). We treat \(2\le\ell\le k-2\), apart from \((k,\ell)=(4,2)\). When \(s<\ell\), we start with an \((\ell+1)\)-graph of prescribed codegree density and add \(s-1\) private vertices to each edge. When \(s\ge\ell\), we instead attach a hypergraph to each transversal of several disjoint copies of a fixed \(k\)-graph. In this range an upper bound comes from the link of a transversal, and a matching construction is obtained by prescribing links on one side of a bipartition.
\subsection{Expansion from codegree}\label{ch:majority}

The expansion starts from a two-part \((\ell+1)\)-graph with \(F\)-free links across the parts. Its codegree controls the minimum \(\ell\)-degree after expansion. For the upper bound we use an auxiliary \((\ell+1)\)-graph of sets with many extensions, and for the lower bound we take cliques in the two-part construction.

\begin{theorem}
For \(k/2<\ell\le k-2\), every real number field \(K\subset\R\) is generated by singleton minimum-\(\ell\)-degree densities tending to one.
\end{theorem}

\begin{proof}
Put \(h=\ell+1\). Since \(s<\ell\) and \(s\ge2\), we have \(\ell\ge3\). By Theorem~\ref{thm:ordinary-fields}, choose a single \(\ell\)-graph \(F\) whose ordinary density \(\beta>2/3\) generates \(K\). We first use \(F\) to construct \(h\)-graphs of a prescribed codegree density.
Partition the vertices into two independent sets \(A,B\). For each \(v\in B\) require the \(\ell\)-graph \(L_H(v)[A]\) to be \(F\)-free, and impose the same condition with \(A\) and \(B\) interchanged. Include all crossing \(h\)-sets with fewer than \(\ell\) vertices in either part.

Suppose that \(A\) is the larger part. Count the edges with \(\ell\) vertices in \(A\) and one in \(B\), first by their vertex in \(B\) and then by their \(\ell\)-set in \(A\). The average codegree of these \(\ell\)-sets is at most \((\beta+o(1))|B|\), so the minimum codegree ratio is at most \(\beta/2+o(1)\).
For the lower bound, balance \(A,B\) and choose an independently relabelled extremal \(F\)-free link for every vertex of the opposite part. By concentration every \(\ell\)-set in a single part has at least \((\beta/2-o(1))n\) neighbours. For a mixed \(\ell\)-set, at least one part supplies all its vertices outside the set as neighbours, since \(\ell\ge3\).

We have thus obtained a two-colour class of codegree density \(\beta/2>1/3\). Theorem~\ref{weakloc:thm:loc} gives a finite forbidden family of this density. Theorem~\ref{cg:thm:conversion} then gives a single forbidden \(h\)-graph \(P\) of codegree density \(\alpha=1-(q+2/(2-\beta))^{-1}\) outside the joined class. In the balanced examples each of the \(q\) added classes has proportion \(\mu=1-\alpha\), and each of the two core classes has proportion \(\mu/(2-\beta)\), up to \(o(1)\).

We can now form a \(k\)-graph from \(P\) by adding \(s-1\) new vertices to each edge, using disjoint sets of new vertices for distinct edges. We show that this expansion has minimum \(\ell\)-degree density \(\alpha^s\).
For the upper bound, let \(H\) be an \(n\)-vertex \(k\)-graph. Form an \(h\)-graph \(J\) on \(V(H)\) by taking as edges those \(h\)-sets \(T\) with \(d_H(T)\ge C n^{s-2}\), for a sufficiently large constant \(C\). For any \(\ell\)-set \(S\), at most \(O(n^{s-1})\) edges through \(S\) use a vertex outside \(S\cup N_J(S)\). All other edges through \(S\) extend it by an \(s\)-set in \(N_J(S)\). Hence, if the minimum \(\ell\)-degree ratio of \(H\) exceeds \(\alpha^s\) by a fixed margin, then the codegree ratio of \(J\) exceeds \(\alpha\) for large \(n\). Embed \(P\) in \(J\), and extend its edges greedily to edges of \(H\) using distinct new vertices.

For the lower bound, let the edges of the \(k\)-graph be the \(k\)-cliques in one of the balanced \(P\)-free \(h\)-graphs above. A copy of the expansion would give a copy of \(P\) in this \(h\)-graph. For an \(\ell\)-set in one class, we may take all \(s<\ell\) additional vertices from its neighbourhood in the \(h\)-graph. The proportion of missing extensions is therefore at most \(1-\alpha^s+o(1)\). A mixed \(\ell\)-set with at most \(\ell-2\) vertices in every class loses \(O_k(\mu^2)\), since a missing edge requires two additional vertices in the same class. Otherwise, exclude the class containing \(\ell-1\) of its vertices. A second monochromatic \(\ell\)-set can arise only when \(s=\ell-1\) and all the additional vertices lie in the other occupied class. All edges crossing distinct components are present. If the two classes belong to the same core, the first excluded class has proportion \(\mu/(2-\beta)<\mu\), leaving enough slack for the additional loss \(O(\mu^s)\). This proves the lower bound. By Lemma~\ref{lem:field-transfers}, \(\alpha^s\) generates \(K\) for all sufficiently large \(q\).
\end{proof}

\subsection{Prescribing links}\label{ch:link}

Here the choice of transversal size determines which result we need in its link. We take \(\ell-1\) vertices when \(\ell\ge3\), so that the link bound is a minimum vertex-degree bound. When \(\ell=2\), we take both vertices of a pair and use ordinary density in its link.

\begin{theorem}
If \(3\le\ell\le s\), or \(\ell=2,\ s\ge3\), every real number field \(K\subset\R\) is generated by singleton minimum-\(\ell\)-degree densities tending to one.
\end{theorem}

\begin{proof}
Let \(d\ge3\). By Section~\ref{sec:ordinary-single} there is a \(d\)-graph \(F\) of ordinary density \(\beta=1-u^{-(d-1)}\), where \(u=q+u_0\), \(\Q(u_0)=K\), and \(q\) can be any sufficiently large integer. By deleting vertices of small degree and then sampling, we obtain \(F\)-free examples with near-extremal minimum vertex degree at all sufficiently large orders. The \(t=0\) suspension \(\mathcal T_0(F)\) has density \(1-(u+1)^{-(d-1)}\), and deleting any independent set from it leaves a copy of \(F\)~\cite{ADLY}. Similarly, suspending a single \(k\)-edge \(r-1\) times yields a \(k\)-graph \(G_r\) of minimum \(\ell\)-degree density \(1-r^{-s}\) and weak chromatic number at least \(r+1\). This proves the assertion for \(K=\Q\). Henceforth assume \(K\ne\Q\),
so \(u_0\) is irrational.

Take \(p\le\ell\) disjoint copies of \(G_r\). For each transversal \(T\) consisting of one vertex from each copy, take a fresh copy of a \((k-p)\)-graph \(P\) and add the edges \(T\cup e\) for all its edges \(e\).
Suppose a host has minimum \(\ell\)-degree ratio above \(\alpha>1-r^{-s}\), where \(\alpha\) is the minimum \((\ell-p)\)-degree density of \(P\) (its ordinary density when \(p=\ell\)). First embed the copies of \(G_r\). For each transversal \(T\), the link \(L_H(T)\) then contains the required copy of \(P\), with previously used vertices excluded.
For the lower bound, partition \(V(H)=A\sqcup B\) so that \(H[B]\) is weakly \(r\)-colourable and \(L_H(T)\) is \(P\)-free for every \(p\)-set \(T\subset A\). Each copy of \(G_r\) must meet \(A\). Choosing one such vertex from every copy would give a transversal in \(A\) whose link contains \(P\), which is impossible.

When \(3\le\ell\le s\), we choose \(p=\ell-1\) so that the link bound uses minimum vertex degree. Take \(d=s+1,\ P=\mathcal T_0(F),\ m=(u+1)^{-1},\ r=\lceil u\rceil\). Since \(q\) is an integer, \(r-u=\lceil u_0\rceil-u_0\) is a fixed
number in \((0,1)\). Thus \(u<r<u+1\), and the upper bound is
\(\alpha=1-m^s\).
Give \(A\) and \(B\) proportions \(m\) and \(um\), respectively. Make \(A\) independent and split \(B\) into \(r\) independent classes of equal size. For every \(p\)-set \(T\subset A\), independently choose a relabelled extremal \(F\)-free hypergraph on \(B\) and use it for \(L_H(T)[B]\). This restriction concerns exactly the edges whose intersection with \(A\) is \(T\). The whole link of \(T\) avoids \(\mathcal T_0(F)\), since its vertices in \(A\) form an independent set and the remaining link is \(F\)-free.

We count the missing extensions of an \(\ell\)-set by its intersection with \(A\). A set contained in \(A\) loses a proportion \(m^s\), and one with \(\ell-1\) vertices in \(A\) loses at most \((1-\beta)(um)^s=m^s\). With \(i\le\ell-2\) vertices in \(A\), the link restrictions cost at most \(\binom{s}{j}u^{-j}m^s+o(1)\), where \(j=\ell-1-i\). This follows by concentration over the independent choices indexed by the additional \(j\) vertices in \(A\), with ranges \(O(n^{s-j})\). An \(\ell\)-set in a single class of \(B\) also loses \((um/r)^s\). Here \(j\ge2\), and the relative slack \(1-(u/r)^s=\Omega(1/u)\) absorbs the loss from the links. This gives equality at \(\alpha\), which generates \(K\) by Lemma~\ref{lem:field-transfers}.

When \(\ell=2\) and \(s\ge3\), we use ordinary density in the links of pairs. Take \(d=s,\ p=2,\ P=\mathcal T_0(F)\), of density \(\alpha=1-\Delta\) with \(\Delta=(u+1)^{-(s-1)}\). Put \(a=(u+1)^{-1}, b=1-a,\ \Lambda=\Delta^{-1/s}\). For each large integer \(r\), take the least shift \(q\) with \(\Lambda>r\). The preceding shift gives \(\Lambda-r\le(s-1)u^{-1/s}/s\), hence \(r-b\Lambda\ge u^{-1/s}/(2s)\); in particular \(1-r^{-s}<\alpha\) and \((b/r)^s\le\Delta(1-\Omega_s(1/u))\).
Take \(A,B\) of proportions \(a,b\) and prescribe an ordinary \(F\)-free \(s\)-graph in the link of each pair in \(A\), as above. A pair in \(A\) loses a proportion \(a^s+(1-\beta)b^s=\Delta\). For a pair meeting both parts the loss is at most \(s\Delta/(u+1)\). For a pair in \(B\), independent relabellings bound the loss from the prescribed links by \(\binom{s}{2}\Delta/[u(u+1)]\). If both vertices are in the same independent class of \(B\), the slack in \((b/r)^s\) absorbs this last term. Thus the density is exactly \(\alpha\), and it generates \(K\) along the chosen shifts by Lemma~\ref{lem:field-transfers}.
\end{proof}

\section{Minimum pair degree in 4-graphs}\label{ch:quartic}

The construction of Section~\ref{ch:link} for \(\ell=2\) prescribes an
ordinary density in the link of a pair. When \(k=4\), such a link is a
graph, whereas the prescribed-field construction of
Section~\ref{sec:ordinary-cases} requires uniformity at least three.
We instead join a core to \(q=4^d\) independent classes indexed by
the points of \(\mathbb F_4^d\). Two distinct labels determine a unique
line, and we impose a parity restriction on the four classes of each
line. We use this incidence property twice, first to count the missing
extensions of a pair and then to rule out competing labels for a row
of the product.

\begin{theorem}[Exact padded conversion]\label{quartic:conv}
Let \(\D\) be a hereditary class of weakly \(p\)-colourable \(4\)-graphs
of minimum \(2\)-degree density \(\beta<1\). Suppose there is a finite
family \(\F\) with \(\pi_2(\F)=\beta\) such that every member of
\(\D\) is \(\F\)-free. For every sufficiently large power \(q=4^d\),
there is a single \(4\)-graph of minimum \(2\)-degree density
\(1-(q+(1-\beta)^{-1/2})^{-2}\).
\end{theorem}

The parameter in this formula is \((1-\beta)^{-1/2}\). Even if \(\beta\)
generates \(K\), this parameter may lie outside \(K\). For example, at
\(\beta=1/2\) it equals \(\sqrt2\). We shall therefore choose the core so that
\(1-\beta=y^2\) for a positive \(y\) with \(\Q(y)=K\).
Lemma~\ref{lem:field-transfers}, applied to \(q+1/y\) with
\(R(X)=1-X^{-2}\), then gives the following result.

\begin{theorem}[Prescribed field in the middle case]\label{quartic:field}
Every real number field \(K\subset\R\) is generated by singleton minimum \(2\)-degree densities of \(4\)-graphs arbitrarily close to one.
\end{theorem}

\subsection{The extremal class}\label{sec:quartic-class}
Put \(c=(1-\beta)^{-1/2}\), \(U=q+c\) and \(\alpha=1-U^{-2}\).
For \(q=4^d\), fix a bit \(s_B\) for each affine line \(B\) in
\(\mathbb F_4^d\), and write \(\mathbf s=(s_B)_B\). A member of
\(\mathcal L_q(\mathbf s)\) has a core in \(\D\) and \(q\) independent
classes \(V_a\) indexed by \(a\in\mathbb F_4^d\). Give each vertex
outside the core a bit. A \(4\)-set meeting the four classes on a line
\(B\) is forbidden when the sum of its bits equals \(s_B\) modulo two.
All other crossing edges are unrestricted. The following density
calculation holds for every fixed assignment \(\mathbf s\).

\begin{localclaim}
The class \(\mathcal L_q(\mathbf s)\) has extremal minimum pair-degree density
\(\alpha\). There is a finite witness family \(\F_q\) of this density
whose members admit no partition into \(q\) independent classes and
an \(\F\)-free core.
\end{localclaim}

\begin{proof}
For the upper bound, apply Lemma~\ref{lem:family-suspension} with
\((k,\ell)=(4,2)\) a total of \(q\) times. It gives a finite family
\(\F_q\) of density \(\alpha\) whose members lie outside joins of \(q\)
pads and an \(\F\)-free core. Indeed, at one step take two disjoint root members. At each pair of root vertices attach a member, and at each vertex of that member attach another, joining it to the root pair and that vertex. An independent set misses a root, or for a selected root pair in it either misses their attachment or forces its second attachment to survive. For input threshold \(\gamma\), set \(\theta=(1+(1-\gamma)^{-1/2})^{-1}\). Above \(1-\theta^2\), each subset of order \((1-\theta)n-O(1)\) has induced ratio above \(\gamma\) after any bounded deletions. In every root-pair link the density margin above \(1-\theta^2\) gives more than \((1-\theta)n\) vertices of graph degree at least \((1-\theta+\epsilon')n\) for small \(\epsilon'>0\) (otherwise bound its density by \(1-\theta^2+\theta\epsilon'+o(1)\)). This embeds both kinds of attachments. Balancing an independent pad of mass \(\theta\) and the core gives the matching threshold.

To attain the bound, give the core proportion \(c/U\) and each \(V_a\) proportion \(1/U\), divided equally between the two bit values. A pair in one component can lose extensions only within that component. A pair in distinct classes \(V_a,V_b\) lies on a unique line of class labels. Of the extensions using one vertex from each of the other two classes on the line, exactly half have forbidden parity, a loss of \(U^{-2}+o(1)\). A pair meeting both the core and a class \(V_a\) loses no extensions. Thus the overall normalized pair-degree deficit is \(U^{-2}+o(1)\).
\end{proof}

\subsection{A single forbidden hypergraph}\label{sec:quartic-single}
Put \(\mu=U^{-2}+U^{-3}/10\) and \(\lambda=\sqrt\mu\). Take the categorical product \(P\) of all \(m\)-vertex \(4\)-graphs that have minimum pair-degree ratio at least \(1-\mu\) and contain a member of \(\F_q\). We choose \(m\) sufficiently large and divisible by the clique order specified below. The witness in each factor will give the upper bound. For the lower bound we exploit the line signs to show that a colouring of the product would give an impossible partition of one factor.

For the lower bound, suppose a row of one factor could receive two
distinct pad labels \(a,b\) at different outer positions. We shall test
both labels against a vertex in a row of a third label \(w\). When the
lines \(aw\) and \(bw\) are distinct, the two tests use different pairs
of completing classes. If both parity restrictions apply, the same pair
of factor vertices loses more completions than its minimum pair degree
allows. The next lemma chooses the line signs so that there are many
such labels \(w\).
\begin{lemma}\label{quartic:sign}
For every sufficiently large \(q=4^d\), there is an assignment
\(\mathbf s\) with the following property. For every
\(\sigma\in\{0,1\}^q\) there is a set \(E_\sigma\) of \(O(1)\) labels
such that, for all distinct \(a,b\notin E_\sigma\) and all \(u,v\in\{0,1\}\),
at least \(q/10\) labels \(w\notin\{a,b\}\) give distinct lines
\(B=L(a,w)\) and \(B'=L(b,w)\) satisfying
\[
 u+\sum_{j\in B\setminus\{a\}}\sigma_j=s_B,\qquad
 v+\sum_{j\in B'\setminus\{b\}}\sigma_j=s_{B'}
 \quad\text{in }\mathbb F_2 .
\]
\end{lemma}
\begin{proof}
Choose the line signs independently and uniformly. Fix \(\sigma\),
\(M\) disjoint pairs and their bits. For each pair discard probes on the
line through the pair and probes whose line through either endpoint
contains an endpoint of another pair. This removes \(O(M)\) probes per
pair, since each line has four points. Each remaining probe satisfies
both equations with probability \(1/4\). Changing one sign affects at
most three probes for that pair, so bounded differences gives probability
\(e^{-\Omega(q)}\) of fewer than \(q/10\) successes. The sets of signs
used by different pairs are disjoint. Thus the probability that all \(M\)
pairs fail is \(e^{-\Omega(Mq)}\). A union bound over at most
\(2^q q^{2M}4^M\) choices gives the assertion for a sufficiently large
constant \(M\). For each \(\sigma\), take \(E_\sigma\) to be the endpoints
of a maximal matching of pairs that fail for some choice of bits.
\end{proof}

The next estimate finds a product edge through a given outer root tuple.
\begin{lemma}[Rooted mixing]\label{quartic:mix}
At any outer root tuple, three subsets each of measure at least \(0.4\) on the leave-one-out product contain a transversal of its product of links.
\end{lemma}
\begin{proof}
In each \(3\)-link the one-position missing mass, including repetitions, is \(O(\mu+1/m)\) by minimum pair degree. Apply the normalization of Section~\ref{ch:product} with the uniform law on the vertices other than the root to get uniform marginals and \(\|M_i-1\|_2=O(\sqrt{\mu+1/m})\). For large \(q\) and then large \(m\) the displayed norm is below \(d_0(3)\). Lemma~\ref{lem:product-correlations} gives positive joint mass.
\end{proof}

\begin{proof}[Proof of Theorem~\ref{quartic:conv}]
For the upper bound, a host above \(\alpha\) contains a bounded witness from \(\F_q\). Keep its vertices and choose the rest of an \(m\)-vertex set at random. Bounded differences and the fixed degree margin give an induced \(4\)-graph with all pair degrees as required. This embeds a factor, so Lemma~\ref{lem:product-upper} gives \(\pi_2(P)\le\alpha\).

For the lower bound, fix one assignment \(\mathbf s\) supplied by
Lemma~\ref{quartic:sign} for this \(q\). Use it for every host in
\(\mathcal L_q(\mathbf s)\). Suppose that \(P\) is a subhypergraph of
some \(H\) in this class.
Restrict the labels and bits of \(H\) to its copy of \(P\), refining the
core colour into \(p\) independent classes.

Let \(f\) be the indicator of an independent set in \(P\) and \(a\) its proportion of the vertices. We claim \(a\le\lambda+O(\lambda^{3/2})\) and \(\tail(f)=O((\lambda-a)_++\lambda^{3/2})\). Fix a vertex in the set and restrict to the product of the \(3\)-uniform links at its coordinates. The factor links have minimum vertex-degree ratio at least \(1-\mu\), so we can apply Lemma~\ref{lem:ordinary-product-tail} with \(s=2\). Transfer via the four-block statistic in Lemma~\ref{cg:lem:tail}, with error \(O(1/(m\sqrt a))\) (use the trivial bound when \(a\le\lambda/2\)). For \(t=p+q\) the mass deficits and total tail are \(O(\sqrt\lambda)\), so Lemma~\ref{lem:one-coordinate} gives a coordinate on which the most frequent fibre colour has error \(O(\sqrt\lambda)\).

Fix this coordinate. In its factor, all but \(O(\sqrt q)\) pads have disjoint good slices \(S_j\) of at least \(0.99\lambda m\) rows on each of which colour \(j\) has probability at least \(0.99\). The exceptional classes are accounted for by the mass deficits, the error of the fibre colouring and the rows discarded in forming the slices. Choose a bit in each row with correct-pad-and-bit probability \(\ge0.495\), and let \(\sigma_j\) be the majority of these bits in the slice.
Assign arbitrary bits to the labels without a good slice. Chebyshev for the exclusion kernel (centred norm \(1/(m-1)\)), followed by Markov and the union bound over slices, leaves a set \(\Gamma\) of outer tuples of measure \(1-O(q/m^2)\) where at most \(0.01\) of each slice's rows have conditional mass below \(0.4\).

Remove the labels in \(E_\sigma\) given by Lemma~\ref{quartic:sign}
for this vector. Suppose a fixed row \(x\) has two remaining good pad labels \(a,b\), with bits \(u,v\) at outer tuples \(z,z'\in\Gamma\). We choose a probe \(w\) with good slice and both completion profiles good (each bad pad rules out \(O(1)\) probes). Choose \(r\ne x\) in \(S_w\) of majority bit passing the tests at both roots, possible among at least \(0.48|S_w|\) rows. In the two other slices \(S_e,S_f\) of line \(aw\), the parity predicted by the majority bits occurs on at least half the pairs: its fraction is \((1+(2h_e-1)(2h_f-1))/2\) for majority proportions \(h_e,h_f\). At least \(0.47|S_e||S_f|\) survive the tests at \(z\) and avoid repetitions. Each is a missing completion of \(x,r\), since otherwise Lemma~\ref{quartic:mix} gives a product edge of forbidden parity. The same holds on \(bw\) at \(z'\); its completion profile is different. Their total \(0.94(0.99\lambda m)^2\) exceeds \(\mu\binom{m-2}{2}\), contradiction.

It follows that, over outer tuples in \(\Gamma\), each row uses at most \(b_q=O(\sqrt q)\) colours counting the core as one. Let \(L_0\ge4\) bound the orders of the witnesses in \(\F\). Take \(L\ge b_qL_0\) of order \(\sqrt q\), with \(L\mid m\). Lemma~\ref{lem:uniform-cliques} applies since \(\mu L^3=O(q^{-1/2})\) and averaging bounds missing vertex degree. All aligned positions may lie in \(\Gamma\). Give each row a colour occurring at least \(L/b_q\) times. An edge in
one pad of the chosen factor would lift at distinct complete outer
positions to an edge inside a pad of the product. Likewise, a member of
\(\F\) in the core of the factor would lift into the core of \(H\), which
belongs to \(\D\) and is \(\F\)-free by hypothesis. These labels therefore partition the chosen factor into \(q\) independent classes and an \(\F\)-free core, contrary to the member of \(\F_q\) that it contains. For a sole factor the witness already excludes it. Hence every member of \(\mathcal L_q(\mathbf s)\) is \(P\)-free, and the
balanced examples give the matching lower bound.
\end{proof}

We now choose the core for Theorem~\ref{quartic:conv}. We give a
finite-family construction for general \(k,\ell\). For our application
we shall use its shifted windows at \((k,\ell)=(4,2)\), where it turns a
\(3\)-graph codegree density \(\gamma\) into a core of minimum pair-degree
density \(2\gamma/m^2\). This explicit factor will let us prescribe
\(1-\beta\) as a square.
\begin{theorem}\label{thm:finite-degree-families}
For \(k>\ell\ge2\) and every number field \(K\subset\R\), some finite forbidden family has minimum \(\ell\)-degree density generating \(K\).
\end{theorem}
Put \(s=k-\ell,\ m\ge k,\ c_m=s!/m^s\). Following Lin and Zhou~\cite{LZ}, partition the vertices into \(m\) classes and require every edge containing an \(\ell\)-set in class \(i\) to use one additional vertex from each of the classes in \(W_i=\{i+1,\ldots,i+s\}\pmod m\).

\begin{lemma}[Finite localization]
For fixed \(k,\ell,m\), a positive finitely described class with this rule and ratio \(\gamma>c_m-\eta_0\), for sufficiently small \(\eta_0>0\), has finite exact witnesses.
\end{lemma}
\begin{proof}
Choose \(\epsilon>0\) with \(sm\epsilon<c_m/4\), then \(\eta_0\) so small that \((c_m-\eta_0)/(s!(s/m-\epsilon)^s)>s^{-s}\). Choose a bounded vertex set \(A\) whose induced hypergraph retains the required minimum degree. In any representation of this anchor, a largest vertex class contains an \(\ell\)-set. Its many extensions force the successor class to be large as well. Continuing around the cycle gives an \(\ell\)-set in every class. Fix a colouring of \(A\) and an \(\ell\)-set \(S_i\) of each colour. Let \(J_i=L_H(S_i)[V(H)\setminus A]\) and let \(U_i\) be the union of its edges. If \(|U_i|<(s/m-\epsilon)n\), the displayed inequality gives \(|E(J_i)|/|U_i|^s>s^{-s}+\Omega(1)\); a bounded sample then violates the AM--GM edge bound for \(s\)-partite graphs and rejects the anchor. Otherwise a vertex's membership in the \(U_i\)'s restricts it to their window intersection, certifiable by one edge per incidence. An empty list rules out the colouring. The count \(\sum|U_i|\ge(s-m\epsilon)n\) pins all but \(m\epsilon n+O(1)\) vertices, since the distinct full incidence patterns have size \(s\). Each unpinned \(\ell\)-set has a completion entirely pinned, by the \(sm\epsilon\) bound. The distinct windows allow at most one pure colour, certifiable by that edge and pins. Any two assignments differ at most \(m(m-1)(\ell-1)\) times. Apply Lemma~\ref{lem:bounded-colour}.
\end{proof}

\begin{proof}[Proof of Theorem~\ref{thm:finite-degree-families}]
First suppose \(\ell\ge3\). Take an ordinary \(\ell\)-graph witness family of density \(\beta\) generating \(K\) and close to one, by Section~\ref{sec:ordinary-cases}. For each transversal \(T\) of the classes in \(W_i\), require the \(\ell\)-graph \(L_H(T)[V_i]\) to avoid that family. Include all edges with fewer than \(\ell\) vertices in any one class. Averaging the degrees of the \(\ell\)-sets contained in each class bounds the degree by \((\beta+o(1))s!\min_i\prod_{j\in W_i}x_j\le c_m\beta+o(1)\); all colours are linear by the successor argument. Balance the classes and independently relabel the near-extremal links.
Concentration gives the required degrees for pure \(\ell\)-sets. A mixed
\(\ell\)-set meets at most \(\ell\) classes. Since \(m\ge k=\ell+s\),
we may choose \(s\) distinct unused classes and add one vertex from each.
Every resulting \(k\)-set has fewer than \(\ell\) vertices in each class,
so these extensions alone give normalized degree at least \(s!/m^s-o(1)\).

Thus the density is \(c_m\beta\), and localization gives the result.

\phantomsection\label{ch:loc}A shifted choice of the cyclic windows gives a second construction and covers \(\ell=2\) as well.
For \(s\ge2\) shift the windows to \(W_i'=\{i,\ldots,i+s-1\}\); the localization is unchanged (successors and distinct incidence patterns remain). The pure profile is now \((\ell+1)\) vertices in \(i\) and one in each other window colour. On \(V_i\), form the \((\ell+1)\)-graph whose edges are all \((\ell+1)\)-sets that occur in an edge of this profile. Require this shadow to avoid an \((\ell+1)\)-graph whose codegree density \(\beta\) is close to one and generates \(K\). This is a finite positive prescription: choose one extension per shadow edge, allowing any overlaps. An \(\ell\)-set of least codegree in each shadow gives a ratio at most \((\beta+o(1))s!\prod_{j\in W_i'}x_j\); multiply the window products as before. Balanced classes with extremal shadows, all their extensions and all edges without pure facets attain \(c_m\beta\). This works identically for a finite input family. The case \(s=1\) is the codegree theorem.
\end{proof}

\begin{proof}[Proof of Theorem~\ref{quartic:field}]
It remains to choose a core for which \(1-\beta\) is a square in \(K\) with a square root generating \(K\). Let \(G\subset K_\infty\) be the open set of generating deficits from Theorem~\ref{ch:region}, and set \(c_0=2/m_0^2\) for fixed \(m_0\ge4\). The preceding shifted-window construction  turns each finite-family \(3\)-codegree \(\gamma\) near one into a weakly \(m_0\)-colourable core of minimum \(2\)-degree \(\beta=c_0\gamma\) with finite exact witnesses. Deleting vertices or edges preserves the
coloured window and shadow restrictions, so every member of this core
class is free of each witness as a subhypergraph.

Choose positive rational \(h\) with \(x=(1-h^2)/c_0=a/e\) just below one. For \(ap\equiv1\pmod e\) set \(b=(ap-1)/e\), fix \(g_*\in K\cap G\), and put
\[
 j=\frac{c_0}{2eh},\qquad v=\frac{c_0 g_*-j^2}{2h},\qquad
 s_p=h+\frac jp+\frac v{p^2}.
\]
The deficit \(g_p=p^2(s_p^2-1+c_0b/p)/c_0\) equals \(g_*+2jv/(c_0p)+v^2/(c_0p^2)\). It lies in \(G\) for large \(p\) and has height \(p^{O(1)}\), so the exact prescription applies to \(\gamma=b/p-g_p/p^2\). Its projected core has \(\beta=1-s_p^2\), where \(s_p>0\) for large \(p\). For fixed \(p\), \(s_p\) is an affine
function of \(g_*\) with nonzero rational slope \(c_0/(2hp^2)\), so
\(\Q(s_p)=\Q(g_*)=K\). Theorem~\ref{quartic:conv} gives \(1-(q+1/s_p)^{-2}\), which generates \(K\) for all large admissible \(q\) by Lemma~\ref{lem:field-transfers}.
\end{proof}

\begin{proof}[Proof of Theorem~\ref{thm:main}]
For \(\ell=1\) apply Theorem~\ref{thm:ordinary-fields}, and for
\(\ell=k-1\) apply Theorem~\ref{cg:main}.
For \(2\le\ell\le k-2\) with \((k,\ell)\ne(4,2)\), the two
constructions in Section~\ref{sec:intermediate-cases} cover
\(k-\ell<\ell\) and \(k-\ell\ge\ell\), respectively.
The case \((k,\ell)=(4,2)\) follows from Theorem~\ref{quartic:field}.
In each case there are single forbidden \(k\)-graphs whose densities
generate \(K\) and tend to one from below.
\end{proof}

\section*{Statement on AI}
The authors used AI as an exploratory tool to assist in developing the arguments. The authors
have carefully checked all mathematical proofs and take full responsibility for the content of the
manuscript.

\bibliographystyle{abbrv}
\bibliography{references}

\end{document}